\documentclass[11pt]{amsart} \textwidth=14.5cm \oddsidemargin=1cm
\usepackage{mathrsfs}
\usepackage{amsmath, amscd, amsthm, amssymb, amsfonts, verbatim,subfigure}

\usepackage[all]{xy}
\usepackage[dvips]{graphicx,color}

\title[TCFTs and Calabi-Yau categories]{Topological conformal field theories and  Calabi-Yau categories}
\author{Kevin Costello}
\address{Department of Mathematics \\Imperial College London}
\email{k.costello@imperial.ac.uk}

\date{}

\newcommand{\K}{\mbb K}
\newcommand{\til}{\widetilde}
\newcommand{\mscr}{\mathscr}
\renewcommand{\det}{\operatorname{det}}

\newcommand{\adm}{\mc A}

\newcommand{\choc}{\mscr{D}}
\newcommand{\ocell}{\mc G}
\newcommand{\oc}{\mc M}
\newcommand{\pc}{\br {\mc N}}
\newcommand{\Comp}{\op{Comp}_{\K}}
\newcommand{\Vect}{\op{Vect}^\ast}

\newcommand{\br}{\overline}

\newcommand{\iso}{\cong}
\newcommand{\diag}{\bigtriangleup} 
\newcommand{\C}{\mathbb C}

\newcommand{\Q}{\mbb Q}

\newcommand{\Z}{\mathbb Z}
\newcommand{\defeq}{\overset{\text{def}}{=}}
\newcommand{\into}{\hookrightarrow}

\newcommand{\op}{\operatorname}

\newcommand{\mbb}{\mathbb}

\newcommand{\mc}{\mathcal}
\newcommand{\from}{\leftarrow}
\newcommand{\ip}[1]{\left\langle #1 \right\rangle}
\newcommand{\abs}[1]{\left| #1 \right|}
\newcommand{\cmod}{\overline {\mc M}}
\renewcommand{\mod}{\op{mod}}

\newcommand{\R}{\mbb R}
\renewcommand{\d}{\mathrm{d}}

\newcommand{\tr}{\triangle}

\newtheorem*{udefinition}{Definition}
\newtheorem*{thmA}{Theorem A}

\newtheorem{conjecture}{Conjecture}
\newtheorem*{ucorollary}{Corollary}

\newtheorem{theorem}{Theorem}[subsection]
\newtheorem{thm-def}{Theorem/Definition}[theorem]
\newtheorem{proposition}[theorem]{Proposition}
\newtheorem{definition}[theorem]{Definition}
\newtheorem{lemma}[theorem]{Lemma}

\numberwithin{equation}{subsection}

\begin{document}

\maketitle
\begin{abstract}

This is the first of two papers which construct a purely algebraic counterpart to the theory of Gromov-Witten invariants (at all genera).  These Gromov-Witten type invariants depend on a Calabi-Yau $A_\infty$ category, which plays the role of the target in ordinary Gromov-Witten theory.  When we use an appropriate $A_\infty$ version of the derived category of coherent sheaves on a Calabi-Yau variety, this constructs the $B$ model at all genera.   When the Fukaya category of a compact symplectic manifold $X$ is used, it is shown, under certain assumptions, that the usual Gromov-Witten invariants are recovered. The assumptions are that open-closed Gromov-Witten theory can be constructed for $X$, and that the natural map from the Hochschild homology of the Fukaya category of $X$ to the ordinary homology of $X$ is an isomorphism.

\end{abstract}

\tableofcontents

\section{Introduction}

If $X$ is a Calabi-Yau manifold, Witten \cite{wit1991b} describes two different topological twistings of the non-linear sigma model of maps from a Riemann surface to $X$, which he calls the $A$ and $B$ models.  If $X,X^{\vee}$ are a mirror pair of Calabi-Yau varieties, then the $A$ model on $X$ is equivalent to the $B$ model on $X^{\vee}$, and vice-versa.  

The $A$ model has been mathematically constructed as the theory of Gromov-Witten invariants.  The genus $0$ part of the $B$ model has been constructed by Barannikov-Kontesvich \cite{bar_kon1997} and Barannikov \cite{bar2000,bar1999}.  They construct a Frobenius manifold from the variations of Hodge structure of a Calabi-Yau.  The genus $0$ part of mirror symmetry is then the statement that the genus $0$ part of the Gromov-Witten theory of a Calabi-Yau variety $X$ is equivalent to the theory of Barannikov-Kontsevich on a Calabi-Yau $X^{\vee}$. 

The higher genus $B$ model is more mysterious.  In the physics literature, it is constructed as a kind of quantisation of the Kodaira-Spencer deformation theory of complex structures on a Calabi-Yau \cite{bcov}.  

However, despite the great deal of interest in mirror symmetry since the subject's inception in the early 1990's, there has been no rigorous construction of the higher-genus part of the $B$ model.  One of the aims of this paper is to construct the $B$ model rigorously for the first time, and so provide a mirror partner to the entire theory of Gromov-Witten invariants.  

\vspace{0.3cm}

Kontsevich \cite{kon1995icm} formulated mirror symmetry as an equivalence of $A_\infty$ categories. If $X,X^{\vee}$ are a mirror pair, then Kontsevich conjectures that the Fukaya category of a variety $X$ ($A$ model) is equivalent to the dg category of complexes of coherent sheaves on $X^{\vee}$ ($B$ model).  Kontsevich's homological mirror symmetry conjecture should explain other aspects of mirror symmetry.  In particular, the equivalence of the theory of Gromov-Witten invariants on $X$ with the $B$ model on $X^{\vee}$ should be a corollary of Kontsevich's conjecture.  

Both of the $A_\infty$ categories appearing in Kontsevich's conjecture are of Calabi-Yau type.  This means, roughly, that there is a non-degenerate invariant pairing on the space of morphisms.

This immediately suggests the following picture.  From each Calabi-Yau $A_\infty$ category, one should construct something like the theory of Gromov-Witten invariants.  If the input Calabi-Yau $A_\infty$ cateogry is the Fukaya category of a compact symplectic manifold, then this theory should recover the usual theory of Gromov-Witten invariants.  If the input Calabi-Yau $A_\infty$ category is the category of sheaves on a smooth projective variety, the resulting theory will, by definition, be the $B$ model at all genera.

\vspace{0.2cm}

In this paper, we prove results along these lines. These results are derived from a study of a kind of abstract topological string theory, called a topological conformal field theory (TCFT).  We study open, closed and open-closed TCFTs. 

Closed TCFTs behave like the Gromov-Witten invariants of a projective variety : a closed TCFT can be described as a collection of cochains on moduli space of Riemann surfaces, with values in tensor powers of an auxiliary chain complex (the complex of ``closed states''), and which satisfy certain gluing constraints.  

The main results of this paper are as follows.  Firstly, we show that open TCFTs are the same as Calabi-Yau $A_\infty$ categories.  Thus to each Calabi-Yau variety we have two open TCFTs : that associated to the Fukaya category ($A$ model) and that coming from coherent sheaves ($B$ model). 

Then, we show how one can associate to each open TCFT an open-closed TCFT, and in particular a closed TCFT.  This is a formal, categorical construction.  We observe that to each open TCFT one can associate the homotopy universal open-closed TCFT (this is an example of a homotopy Kan extension).   Then we calculate the homology of the space of closed states of this universal closed TCFT : it is the Hochschild homology of the $A_\infty$ category associated to the open TCFT. 

Also, we show, under certain assumptions, how to relate the closed TCFT constructed here from the Fukaya category of a compact symplectic manifold to the ordinary Gromov-Witten invariants of the manifold.  
 
\vspace{0.3cm}

These results are proved using a combination of homotopical algebra, and some results about the topology of the moduli spaces of Riemann surfaces.  In particular, the dual version of the ribbon graph decomposition of moduli space \cite{cos_2006} plays an essential role.

Let us now turn to describing these results in more detail. 

\subsection{Topological conformal field theories}

Let $\mc M$ be Segal's category of Riemann surfaces.  The objects
of $\mc M$ are finite sets; for sets $I,J$, a morphism from $I$ to
$J$  is a Riemann surface with $I$ incoming and $J$ outgoing
parameterised boundary components. (We require that there is at
least one incoming boundary on each component). Composition of
morphisms is given by gluing of Riemann surfaces. Disjoint union
of sets and of surfaces gives $\mc M$ the structure of symmetric
monoidal category. According to Segal \cite{seg1988}, a conformal
field theory is a symmetric monoidal functor from this category to the
category of vector spaces.

This definition can be modified in several ways.  For example, we
could look for functors from $\mc M$ to the symmetric monoidal category of
topological spaces, or of spectra.   There is a natural linearised
version of these topological functors, obtained by passing from
the category of topological spaces  to the category of chain
complexes. Let $C_\ast$ be a symmetric monoidal functor from the category of
topological spaces to that of complexes of $\K$ vector spaces,
which computes homology groups. (Here $\K$ is a base field of
characteristic zero). The category $\mc M$ has discrete set of
objects, but the spaces of morphisms are topological spaces.
Applying $C_\ast$ to the topological category $\mc M$ yields a
differential-graded category $C_\ast(\mc M)$. The objects of
$C_\ast(\mc M)$ are, as before finite sets; the morphisms of
$C_\ast(\mc M)$ are defined by
$$
\op{Mor}_{C_\ast(\mc M)} (a,b) = C_\ast( \op{Mor}_{\mc M}(a,b) )
$$
Define
$$
\mscr{C} \defeq C_\ast(\mc M)
$$
$\mscr{C}$, like $\mc M$, is a symmetric monoidal category. The
following definition is due independently to Getzler
\cite{get1994} and Segal \cite{seg1999_2}.
\begin{udefinition}
A topological conformal field theory is a symmetric monoidal functor $F$ from
the differential graded category $\mscr{C}$ to the category of
chain complexes.
\end{udefinition}
What this means is the following.  For any finite set, $F(I)$ is a
chain complex. Since $F$ is a symmetric monoidal functor, there is a map
\begin{equation*}
F(I) \otimes F(J) \to F(I \amalg J) \label{eqn symmetric monoidal functor}
\end{equation*}
Usually these maps are required to be isomorphisms; if this was
the case, the functor $F$ would be called \emph{split}.  We relax
this to the condition that these maps are quasi-isomorphisms; we
say the functor is \emph{h-split} (homologically split). Each
chain $\alpha$ in the moduli space of Riemann surfaces with $I$
labelled incoming and $J$ labelled outgoing boundary components
gives a map
$$
F(\alpha) : F(I) \to F(J)
$$
which is of the same degree as $\alpha$.  This map respects the
differential : $F(\d \alpha) = \d F(\alpha)$, where $F(\alpha)$ is
considered as an element of the chain complex $\op{Hom}(F(I),
F(J))$.  Gluing Riemann surfaces together must correspond to
composition of maps, and disjoint union  corresponds to tensor
product.

We need to twist the definition of TCFT by a local system.  Let
$\det$ be the locally constant sheaf of $\K$ lines on the morphism
spaces of the category $\mc M$ whose fibre at a surface $\Sigma$
is
$$\op{det}(\Sigma) = (\op{det} H^\ast(\Sigma) )[-\chi(\Sigma)]$$
This is situated
in degree $\chi(\Sigma)$. If $\Sigma_1, \Sigma_2$ are two surfaces
with the incoming boundaries of $\Sigma_2$ identified with the
outgoing boundaries of $\Sigma_1$, then there is a natural
isomorphism
$$
\det(\Sigma_2 \circ \Sigma_1 ) \iso \det(\Sigma_2) \otimes \det(\Sigma_1)
$$
This shows that if we take chains with local coefficients,
$C_\ast(\mc M, \det)$, then we still get a category.  Let
$$
\mscr{C}^d = C_\ast(\mc M,\det^d)
$$
where we use the notation $\det^d$ for $\det^{\otimes d}$.
\begin{udefinition}
A $d$-dimensional topological conformal field theory is a symmetric monoidal functor from the category $\mscr{C}^d$ to the category of
complexes.
\end{udefinition}
It turns out that the local system $\det$ is trivial (up to a shift). However,
it is still important to keep track of it, especially when we consider
open-closed conformal field theory; although the local system is trivial, in
the open-closed case it can not be trivialised in a way compatible with the
category structure.  In the closed case, this local system is not so important;
however, it is convenient to use it to keep track of the grading.

One apparent disadvantage of the definition of TCFT is that it seems to depend
on an arbitrary choice, that of a chain model for the category $\mc M$.
However, we show that quasi-isomorphic categories have homotopy equivalent (in
a precise sense) categories of functors, so that up to homotopy there is no
ambiguity.

\subsection{Open and open-closed TCFTs}
Open-closed conformal field theory was first axiomatised by Moore
and Segal \cite{moo2001,seg1999_2}.  A Riemann surface with
open-closed boundary is a Riemann surface $\Sigma$, some of whose
boundary components are parameterised, and labelled as closed
(incoming or outgoing); and with some intervals (the open
boundaries) embedded in the remaining boundary components. These
are also parameterised and labelled as incoming and outgoing.  The
boundary of such a surface is partitioned into three types: the
closed boundaries, the open boundary intervals, and the free
boundaries.  The free boundaries are the complement of the closed
boundaries and the open boundary intervals, and are either circles
or intervals. We require that each connected component of $\Sigma$
has at least one free or incoming closed boundary.

\begin{figure}

\includegraphics{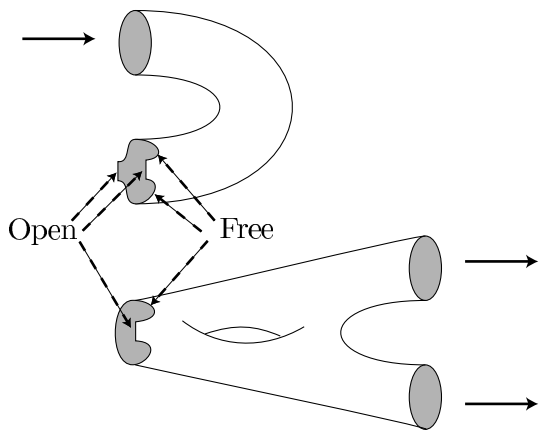}

\caption{A Riemann surface with open-closed boundary.  The open boundaries can
be either incoming or outgoing boundaries, but this is  not illustrated.}
\end{figure}

\begin{figure}
\label{figure open gluing}
\includegraphics{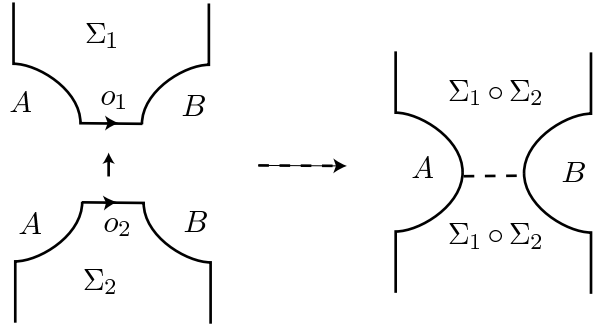}
\caption{ Open gluing, corresponding to composition. $o_1$ on
$\Sigma_1$ is incoming, $o_2$ on $\Sigma_2$ is outgoing, and
$s(o_1) = s(o_2) = A$, $t(o_1) = t(o_2) = B$. Note incoming and
outgoing boundaries are parameterised in the opposite sense.}
\end{figure}

To define an open closed conformal field theory, we need a set
$\Lambda$ of D-branes.  Define a category $\mc M_\Lambda$, whose
objects are pairs $O,C$ of finite sets and maps $s,t : O \to
\Lambda$.  The morphisms in this category are Riemann surfaces
with open-closed boundary, whose free boundaries are labelled by
D-branes.   To each open boundary $o$ of $\Sigma$ is associated an
ordered pair $s(o),t(o)$ of D-branes, where it starts and where it
ends.  Composition is given by gluing of surfaces; we glue all the
outgoing open boundaries of $\Sigma_1$ to the incoming open
boundaries of $\Sigma_2$, and similarly for the closed boundaries,
to get $\Sigma_2 \circ \Sigma_1$.  Open boundaries can only be
glued when their D-brane labels are compatible, as in figure
\ref{figure open gluing}. Disjoint union makes $\mc M_\Lambda$
into a symmetric monoidal category.

Define an open-closed conformal field theory to be a symmetric monoidal functor from $\mc M_\Lambda$ to the category of vector spaces. Let
us assume, for simplicity, that this is split, so that  the
morphisms $F(\alpha) \otimes F(\beta) \to F(\alpha \amalg \beta)$,
for $\alpha, \beta \in \op{Ob} \mc M_{\Lambda}$, are isomorphisms.

Then an open-closed CFT consists of vector spaces $\mc H$, of
closed states; and for each pair of D-branes $\lambda,\lambda'$, a
vector space $\op{Hom} (\lambda,\lambda')$.

Let $\Sigma$ be a Riemann surface with open-closed boundary, each
of whose free boundaries is labelled by a D-brane.  Suppose the
sets of incoming and outgoing closed and open boundaries of
$\Sigma$ are $C_+,C_-,O_+,O_-$ respectively. Then $\Sigma$ must
give a morphism
$$
\mc H( C_+) \otimes \bigotimes_{o \in O_+} \op{Hom}(s(o),t(o)) \to
\mc H( C_-) \otimes \bigotimes_{o \in O_-} \op{Hom}(s(o),t(o))
$$

\begin{figure}

\includegraphics{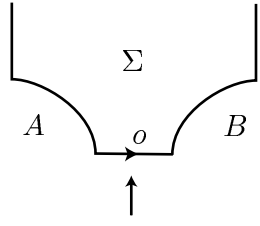}

\caption{$A,B$ are D-branes, labelling free boundaries.   $o$ is
an incoming open boundary with $s(o) = A$, $t(o) = B$.}
\end{figure}
As before, disjoint union of surfaces corresponds to tensor
products of morphisms, and gluing of surfaces -- composition in
the category $\mc M_\Lambda$ -- corresponds to composition of
linear maps.

An open CFT is like this, except  the surfaces have no closed
boundaries, and there is no space of closed states.

The definition of open-closed (or just open) topological  CFT is
obtained from this definition in the same way the definition of
topological CFT is obtained from the definition of CFT.  So we
replace the category $\mc M_\Lambda$ by its associated category of
chains, $C_\ast \mc M_\Lambda$.  We can also take chains with
twisted coefficients; define
$$
\mscr{OC}^d_{\Lambda} = C_\ast(\mc M_{\Lambda},\det^d)
$$
Here $\det$ is a certain local system on the moduli spaces of Riemann surfaces
with open-closed boundary.

 An open-closed TCFT of dimension $d$  is a symmetric monoidal functor from
$\mscr{OC}^d_{\Lambda}$ to complexes, which is h-split, so that the maps
$\Phi(\alpha) \otimes \Phi(\beta) \to \Phi(\alpha \amalg \beta)$ are
quasi-isomorphisms.

 Let $\mscr{O}^d_{\Lambda}$ be the full subcategory whose objects
are purely open; so they are of the form $(C,O)$ where $C =
\emptyset$.  Morphisms in $\mscr{O}^d_{\Lambda}$ are chains on
moduli of surfaces with no closed boundary. An open TCFT is a h-split symmetric monoidal functor from
$\mscr{O}^d_{\Lambda}$ to complexes. 

\subsection{Statement of the main results}
There are functors
$$
i : \mscr{O}^d_{\Lambda} \to \mscr{OC}^d_{\Lambda} \from \mscr{C}^d : j
$$
Let $\Phi$ be an open TCFT, so that $\Phi : \mscr{O}^d_{\Lambda}
\to \Comp$ is a symmetric monoidal functor.  Then we can push forward to get
$i_\ast \Phi : \mscr{OC}^d_{\Lambda} \to \Comp$.  Here $i_\ast$ is
the left adjoint to the pull-back functor
$$i^\ast : \op{Fun}
(\mscr{OC}^d_{\Lambda},\Comp) \to \op{Fun}
(\mscr{O}^d_{\Lambda},\Comp )$$ (here $\Comp$ is the category of
complexes of $\K$ vector spaces).

If we think of a category as like an algebra, then a functor from
a category to complexes is like a (left) module; and we can write
this as
$$
i_\ast \Phi = \mscr{OC}^d_{\Lambda} \otimes_{\mscr{O}^d_{\Lambda}}
\Phi
$$
The functor $i_\ast$ is not exact; it doesn't take
quasi-isomorphisms to quasi-isomorphisms.  Instead, we use the
left derived version
$$
\mbb {L} i_\ast \Phi = \mscr{OC}^d_{\Lambda} \otimes^{\mbb
L}_{\mscr{O}^d_{\Lambda}} \Phi
$$
which is exact.  This is obtained by first replacing $\Phi$ by a
flat resolution, and then applying $i_\ast$.

It turns out that $\mbb{L} i_\ast \Phi$ is an open-closed TCFT
(that is, it is h-split). $\mbb{L} i_\ast \Phi$ is the homotopy
universal open-closed TCFT associated to $\Phi$.

We can pull back along $j$, to get a closed TCFT $j^\ast \mbb{L}
i^\ast \Phi$.  This defines a functor from open to closed TCFTs.
We can think of this functor as taking an open TCFT $\Phi$, and
tensoring with the $\mscr{C}^d - \mscr{O}^d_{\Lambda}$ bimodule,
$\mscr{OC}^d_{\Lambda}$; that is,
$$ j^\ast \mbb{L} i_\ast \Phi = \mscr{OC}^d_{\Lambda}
\otimes^{\mbb L}_{\mscr{O}^d_{\Lambda}} \Phi
$$ considered as a left $\mscr{C}^d$ module.

In this paper the following theorem is proved.
\begin{thmA}
\begin{enumerate}
\item
The category of open TCFTs of dimension $d$, with fixed set of
D-branes $\Lambda$, is homotopy equivalent to the category of
(unital) extended Calabi-Yau $A_\infty$ categories of dimension
$d$, with set of objects $\Lambda$.
\item
For any open TCFT $\Phi$, the homotopy-universal functor
$\mbb{L}i_\ast \Phi : \mscr{OC}^d_{\Lambda} \to \Comp$ is h-split,
and so defines an open-closed TCFT.
\item
Let $HH_\ast(\Phi)$ denote the Hochschild homology of the $A_\infty$
category associated to $\Phi$ by part (1). Then the homology of the
closed states of the open-closed TCFT $\mbb{L} i_\ast \Phi$ is
$HH_\ast(\Phi)$. The homology of the open states is just that of
$\Phi$.

More precisely,  for an object $(O,C) \in \op{Ob}
\mscr{OC}_\Lambda^d$, where $O \in \op{Ob} \mscr{O}_\Lambda^d$, $C
\in \op{Ob} \mscr{C}^d $ (so that $C$ is a finite set), we have
$$
H_\ast ( (\mbb{L} i_\ast \Phi ) (O,C) ) = H_\ast (\Phi(O)) \otimes
HH_\ast(\Phi)^{\otimes C}
$$
In particular, the closed TCFT $j^\ast \mbb{L} i_\ast \Phi$ has
homology
$$
H_\ast ( (j^\ast \mbb{L} i_\ast \Phi ) (C)) = HH_\ast(\Phi)^{\otimes
C}
$$
\end{enumerate}
\end{thmA}

\begin{ucorollary}
The homology of moduli spaces acts on the Hochschild homology of any
Calabi-Yau $A_\infty$ category $\mc D$.   That is there are
operations
$$
H_\ast(\mc M(I,J), \det^d) \otimes HH_\ast(\mc D)^{\otimes I} \to
HH_\ast(\mc D)^{\otimes J}
$$
\end{ucorollary}

Part (1) can be viewed as a categorification of the ribbon graph
decomposition of moduli spaces.   The proof relies on the dual
version of the ribbon graph decomposition proved by the author in
\cite{cos_2004, cos_2006}.   The statement that the categories are
homotopy equivalent has a precise meaning. It means that there are
functors from  open TCFTs to extended CY $A_\infty$ categories, and
from extended CY $A_\infty$ categories to open TCFTs, which are
inverse to each other, up to quasi-isomorphism. A Calabi-Yau
category is the categorical generalisation of a Frobenius algebra.
In a Calabi-Yau $A_\infty$ category, the product is only associative
up to homotopy, and there is a cyclic symmetry condition on the
inner product with the higher multiplications $m_n$. The adjective
``extended'' refers to a small technical generalisation of this
definition which will be explained in section \ref{section proof}.

The homotopy universal closed TCFT $\mbb{L} i_\ast \Phi$ has the property that
for every open-closed TCFT $\Psi$, with a map $\Phi \to i^\ast \Psi$ in an
appropriate homotopy category of TCFTs, there is a map $\mbb{L} i_\ast \Phi \to
\Psi$.  Here $i^\ast \Psi$ is the open TCFT associated to $\Psi$ by forgetting
the closed part; the fact that $\mbb{L} i_\ast \Phi \to \Psi$ is a map of
open-closed TCFTs means that the diagrams
$$
\xymatrix{ \mscr{OC}^d_{\Lambda}(\alpha,\beta) \otimes \mbb{L}
i_\ast \Phi  (\alpha) \ar[r] \ar[d]  &
\mscr{OC}^d_{\Lambda}(\alpha,\beta) \otimes \Psi  (\alpha) \ar[d] \\
\mbb{L} i_\ast \Phi(\beta) \ar[r] & \Psi(\beta)
 }
$$
commute, for all objects $\alpha,\beta$ of
$\mscr{OC}^d_{\Lambda}$.

Passing to homology of the closed states, we see that in
particular, for all finite sets $I,J$, the diagram
$$
\xymatrix{ H_\ast(\mc M(I,J), \det^d) \otimes
HH_\ast(\Phi)^{\otimes I} \ar[r] \ar[d] &  H_\ast(\mc
M(I,J),\det^d) \otimes
H_\ast(\Psi)^{\otimes I} \ar[d] \\
HH_\ast(\Phi)^{\otimes J} \ar[r] & H_\ast(\Psi)^{\otimes J}}
$$
commutes.  Here, $HH_\ast(\Phi)$ refers to the Hochschild homology
of the $A_\infty$ category associated to $\Phi$, under the
correspondence between $A_\infty$ categories and open TCFTs.
$H_\ast(\Psi)$ means $H_\ast(\Psi(1))$, so that $H_\ast(\Psi(I)) =
H_\ast(\Psi)^{\otimes I}$.  $\mc M(I,J)$ is the moduli space of
Riemann surfaces with $I$ incoming and $J$ outgoing boundaries.

One could hope that part (3) of this result should give a natural
algebraic characterisation of the category of chains on moduli
spaces of curves,  as morphisms in some homotopy category between
the functors which assign to a Calabi-Yau $A_\infty$ category the
tensor powers of its Hochschild chains.

\subsection{Relation with Deligne's conjecture}

Theorem A implies a higher genus generalisation of Deligne's
Hochschild cochains conjecture. Deligne conjectured that there is a
homotopy action of the chain operad of the little discs operad on
the Hochschild cochain complex of an algebra.  This has now been
proved by several authors, \cite{tam1998, kon1999, kon_soi2000,
mcc_smi2002}.

A variant of Deligne's conjecture states that the framed little
discs operad acts on the Hochschild cochains of a Frobenius algebra.
This has been proved by Kaufmann \cite{kau2004} and
Tradler-Zeinalian \cite{tra_zei2004}.

The framed little discs operad is the operad of genus zero Riemann
surfaces with boundary. What is shown here is that there is a
homotopy action of chains on all-genus moduli spaces of Riemann
surfaces on the Hochschild chains of a Calabi-Yau $A_\infty$
category, or in particular, of a Frobenius algebra. Restricting to
Riemann surfaces of genus zero with precisely one input, we find a
homotopy co-action of the framed little discs operad on the
Hochschild chain complex.  The Hochschild cochain complex of a
Calabi-Yau $A_\infty$ category is dual to the Hochschild chain
complex.  Therefore we can dualise the coaction on Hochschild chains
to find that the Hochschild cochain complex has a homotopy action of
the framed little discs operad, recovering the result of Kaufman and
Tradler-Zeinalian.

It is not difficult to check that the coproduct on Hochschild homology constructed here, which comes from the class of a point in the moduli space of genus $0$ surfaces with one incoming and two outgoing boundaries, coincides with the dual of the standard cup product on Hochschild cohomology.   Also, the operator on Hochschild homology which comes from the generator of $H_1$ of the moduli space of annuli with one incoming and one outgoing boundary  coincides with the B operator of Connes.  

On the other hand, very few of other operations we construct on Hochschild homology admit such a simple description.  In particular, the \emph{product} we construct on Hochschild homology, which can be described explicitly, seems not to have been considered in the literature before.  

\vspace{0.3cm}

Statements close to the higher-genus analog of Deligne's conjecture proved here have been
conjectured by Kontsevich as far back as 1994 \cite{kon1995icm}, and
have also been conjectured by Segal, Getzler, Kapustin and Rozansky
\cite{kap_roz2004}....  A different approach to constructing an
action of chains on moduli space on Hochschild chains has been
outlined by Kontsevich \cite{kon2003} using the standard ribbon graph
decomposition of moduli spaces of curves, in a lecture at the Hodge
centenary conference in 2003.

\subsection{Relations to the work of Moore-Segal and Lazariou}
Moore and Segal \cite{moo2001} and independently Lazaroiu
\cite{laz2000} have obtained descriptions of open-closed topological
field theories.  Topological field theory (TFT) is a greatly
simplified version of the topological conformal field theory
considered in this paper. Instead of taking the singular chains on
moduli spaces, in TFT we only use $H_0$, or equivalently only
consider topological surfaces (with no conformal structure).

These authors show that an open-closed TFT consists of a not
necessarily commutative Frobenius algebra $A$, a commutative
Frobenius algebra $B$, with a homomorphism from
\begin{equation}
\iota_\ast : B \to Z(A)\label{eqn ms map}
\end{equation}
(where $Z(A)$ is the centre of $A$), satisfying an additional
constraint, called the Cardy condition.

Their result is closely related to ours.   We show that an open
topological conformal field theory, with one D-brane, is the same as
a Frobenius $A_\infty$ algebra. This is obviously the derived, or
homotopy, version of their result.

In our situation, the closed states are not just a Frobenius
algebra. They have a much richer structure coming from the topology
of moduli spaces.  Also, the inner product on the space of closed
states may be degenerate, even on homology.  This is because in this
paper we need the restriction that all of our Riemann surfaces have
at least one incoming boundary, whereas in Moore and Segal's work
this is not imposed.

Suppose $(A,V)$ is an open-closed TCFT, for simplicity with one
D-brane.  Then $A$ is an $A_\infty$ Frobenius algebra.  Then the map
$ HH_\ast(A) \to H_\ast(V)$ we construct is an analog of the map
\ref{eqn ms map}. As, if we dualise we get a map
$$
H_\ast(V)^{\vee} \to HH^\ast(A) = HH_\ast(A)^{\vee}
$$
This map is compatible with the operations coming from the homology
of moduli spaces of curves, so in particular, it is a ring
homomorphism. If $A$ purely of degree zero, and all higher products
vanish, then $HH^0(A)$ is the centre of $A$.  We can view
$HH^\ast(A)$ as a derived analog of the centre, and this map
corresponds to the one constructed by Moore-Segal and Lazaroiu.

The Cardy condition automatically holds in our setting (as it comes
from one of the diagrams of open-closed TFT).  However it holds in a
slightly different form to that used by Moore-Segal and Lazaroiu.
For us, the Cardy condition is expressed in terms of the relation
between the inner product on $A$ and a natural inner product on
$HH_\ast(A)$ (and in particular on $HH_0(A) = A/[A,A]$).  For
Moore-Segal and Lazaroiu, the Cardy condition expresses the relation
between an inner product on $B$, which maps to $HH^0(A)$, and that
on $A$.  If the inner product on $HH_\ast(A)$ was non-degenerate,
then the dual inner product on $HH^\ast(A)$ would satisfy
Moore-Segal's and Lazaroiu's form of the Cardy condition.  However,
the inner product on $HH_\ast(A)$ is often degenerate.

\subsection{The non-unital version of the result}

There is a variant of the main result that deals with non-unital
Calabi-Yau $A_\infty$ categories.  This version is perhaps more
suited to applications, as non-unital Calabi-Yau $A_\infty$
categories are easier to construct and have a better-behaved
deformation theory.

As the proof of the non-unital version is essentially the same, I
will just indicate how the statement differs.  The identity
morphisms on an object of a Calabi-Yau $A_\infty$ category
corresponds to the disc with a single open boundary, and free
boundary labelled by a D-brane.  Therefore, if we want to work with
non-unital Calabi-Yau $A_\infty$ categories, we must remove these
morphisms from the category $\mscr O_\Lambda^d$.  Thus let $\til
{\mscr O}_\Lambda^d \subset \mscr{O}_\Lambda^d$ be the subcategory
such that each connected component of the morphism surfaces is not a
disc with $\le 1$  open boundary.

It turns out we have to perform a similar modification in the closed
case.  That is, let $\til{\mscr{C}}^d \subset \mscr{C}^d$ be the
subcategory such that each connected component of the morphism
surface is not a disc with one incoming closed boundary.  We also
have $\til{\mscr{OC}}^d_\Lambda$, where we disallow  surfaces with a
connected component which is a disc with either $\le 1$ open or $1$
closed boundaries, or an annulus with one closed, one free, and no
open boundaries.

Then the analog of theorem A holds. That is, open TCFTs using $\til
{\mscr O}_\Lambda^d$ are homotopy equivalent to non-unital extended
Calabi-Yau $A_\infty$ categories of dimension $d$.  For each such
variant open TCFT, there is a homotopy-universal open closed (using
$\til{\mscr{OC}}^d_\Lambda$).  The homology of the closed states of
this is the Hochschild homology of the Calabi-Yau $A_\infty$
category associated to the open TCFT.  Here we have to be careful
with the definition of Hochschild homology; for a non-unital
category, the correct definition is to formally augment the category
by adding on unit morphisms, and then quotient out by the subcomplex
of the Hochschild chain complex spanned by these identity morphisms
(considered as Hochschild zero chains).  It is this version of
Hochschild homology we find.

\subsection{Outline of the proof of theorem A}
There are two parts to the proof of the main theorem : a homological
algebraic part, and a geometrical part.

The algebraic part consists of constructing some very general
homotopy theory for functors from differential graded symmetric
monoidal categories.  If $A$ is a differential graded symmetric
monoidal (dgsm) category, we consider a dg symmetric monoidal
functor $F : A \to \op{Comp}_\K$ as a left $A$ module. We define the
notion of tensor product and homotopy tensor product of an $A-B$
bimodule with a $B-C$ bimodule.   The main technical point here is
the result that in certain situations flat resolutions of modules
exist.  These results allow us to show that if $A \to B$ is a
quasi-isomorphism of dgsm categories, then the categories of $A$
modules and $B$ modules are homotopy equivalent.

The geometric part of the proof amounts to giving an explicit
generators-and-relations description for a category quasi-isomorphic
to the category $\mc O_\Lambda$, and for  the category $\mc
{OC}_\Lambda$, considered as a right $\mc O_\Lambda$ module.  These
explicit models are derived from certain cell complexes weakly
equivalent to moduli spaces of Riemann surfaces, constructed in
\cite{cos_2004, cos_2006}.  The cell complexes are compatible with
the open gluing maps, but \emph{not} with the closed gluing maps. At
no point do we construct a cellular model for Segal's category.  I
believe that such a model cannot be constructed using the standard
ribbon graph decomposition.

Let me describe briefly how these cellular models for moduli space
are constructed.  A detailed account is contained in
\cite{cos_2006}. Consider the moduli space $\mc N_{g,h,r,s}$ of
Riemann surfaces of genus $g$, with $h$ boundary components, $r$
marked points on the boundary, and $s$ marked points in the
interior.  The boundary marked points will play the role of open
boundary components, and the marked points in the interior (after we
add the data corresponding to the parameterisation) will play the
role of closed boundary components. We use a partial
compactification $\pc_{g,h,r,s}$ into an orbifold with corners,
whose interior is $\mc N_{g,h,r,s}$. This partial compactification
is modular; it parameterises Riemann surfaces possibly with nodes on
the boundary. These nodes appear when we glue together two boundary
marked points. This operation is homotopic to the operation of
gluing two parameterised intervals on the boundary of a surface
together, which gives the composition in the category of Riemann
surfaces with open boundaries.

Inside $\pc_{g,h,r,s}$ is an orbi-cell complex $D_{g,h,r,s}$, which
parameterises Riemann surfaces glued together from discs, each of
which has at most one internal marked point.  This cell complex is
compatible with the open gluing maps; if we take a Riemann surface
built from discs, and glue two of the marked points, the surface is
still built from discs.

If we pass to cellular chains, and restrict to the surfaces with no
internal marked points, we can construct a chain model for the
category $\mscr O_\Lambda$.  It turns out the generators are discs,
and there are only some very simple relations.  The compactified
moduli space of marked points on the boundary of a disc is a
Stasheff polytope.  From this we deduce that open TCFTs are homotopy
equivalent to Calabi-Yau $A_\infty$ categories.

From considering the moduli spaces $\mc D_{g,h,r,s}$ where $s \ge
0$, we can find a model for $\mscr{OC}_\Lambda$ as a right
$\mscr{O}_\Lambda$ module.  This again has a very simple generators
and relations description.  The generators are annuli, one of whose
boundaries is a closed (parameterised) boundary, and the other has
some open marked points on it. (We get annuli from discs with a
single internal marked point, by fattening this marked point into a
(parameterised) closed boundary.  Up to homotopy there is an $S^1$
of ways of doing this).  There is only one relation, which tells us
about forgetting marked points on the boundary of the annulus.

This model allows us to compute the homology of $\mscr {OC}_\Lambda
\otimes_{\mscr{O}_\Lambda}^{\mbb L} F$, for any open TCFT $F$.
  We find this is
the Hochschild homology of the $A_\infty$ category associated to
$F$. This turns out to follow from simple facts about the geometry
of the compactified moduli space of marked points on the boundary of
an annulus.

\section{Examples and applications}

A Calabi-Yau category is the categorical generalisation of a
Frobenius algebra. A CY category $\mc C$ of dimension $d$ (over our
base field $\K$) is a linear category with a trace map
$$
\op{Tr}_A  : \op{Hom}(A,A) \to \K[d]
$$
for each object $A$ of $\mc C$. The associated pairing
$$
\ip{\quad}_{A,B} : \op{Hom}(A,B) \otimes \op{Hom}(B,A) \to \K[d]
$$
given by $\op{Tr} (\alpha \beta)$ is required to be symmetric and
non-degenerate. A Calabi-Yau category with one object is then the same as a
Frobenius algebra.  The grading convention is slightly funny; note that
$\op{Hom}_i(A,B)$ is dual to $\op{Hom}_{ - d - i} ( B,A)$.  This is forced on
us by using homological grading conventions, so the differential is of degree
$-1$.

A Calabi-Yau $A_\infty$ category is an $A_\infty$ category with a trace map as
above, whose associated pairing is symmetric and non-degenerate.  If $\alpha_i
: A_i \to A_{i+1 \op{mod} n}$ are morphisms, then $\ip{m_{n-1}(\alpha_0 \otimes
\ldots \otimes \alpha_{n-2}), \alpha_{n-1} }$ is required to be cyclically
symmetric.

The notion of extended CY $A_\infty$ category is a small technical
generalisation of this definition, and will be explained later.

One special property enjoyed by Calabi-Yau $A_\infty$ categories is a duality
between Hochschild homology and cohomology;
$$
HH_i(D) \iso HH^{d+i} (D)^{\vee}
$$
where $d$ is the dimension of the category.

Our main result implies that the homology of moduli space acts on
the Hochschild homology groups of an Calabi-Yau $A_\infty$ category.
 Next we will discuss in detail what happens for some naturally
arising classes of Calabi-Yau $A_\infty$ categories, associated to a
compact oriented manifold, a smooth projective Calabi-Yau variety,
or a symplectic manifold.

\subsection{String topology}

Let $M$ be a compact, simply connected, oriented manifold.   The cohomology of $M$ has the structure of $C_\infty$ (homotopy commutative) algebra, encoding the rational homotopy type of the manifold.  Hamilton and Lazarev \cite{ham_laz2004} have shown how this enriches naturally to a Frobenius $C_\infty$ algebra, that is a $C_\infty$ algebra with a non-degenerate invariant pairing\footnote{Hamilton and Lazarev's main result is that the deformation theory for Frobenius $C_\infty$ and $C_\infty$ algebras coincide; they deduce the existence of the Frobenius $C_\infty$ structure as an immediate corollary. Note that in the associative world, Frobenius $A_\infty$ and $A_\infty$ algebras have different deformation theory.} . The pairing is simply the Poincar\'e pairing. 

Thus, $H^{-*}(M)$ \footnote{All our complexes are homological, so
we reverse the usual grading} is, in a natural way, a Calabi-Yau $A_\infty$ category with one object. 

Since $H^{-\ast}(M)$ is quasi-isomorphic, as an $A_\infty$ algebra,
to $\Omega^{-\ast}(M)$, a well-known theorem of Adams-Chen implies
that
$$
HH_\ast(H^{-\ast}(M)) = H^{-\ast}(\mc L M)
$$
is the cohomology of the free loop space $\mc L M$ of $M$.

Theorem A now implies that the homology of the moduli spaces of Riemann
surfaces acts on $H^{-\ast}(\mc LM)$.  That is, there are maps
$$
H_\ast(\mc M(I,J),\det^d) \otimes H^{-\ast}(\mc LM)^{\otimes I} \to
H^{-\ast}(\mc LM)^{\otimes J}
$$
compatible with composition and disjoint union.  These operations
should correspond to the higher-genus version of the string topology
operations of Chas-Sullivan \cite{cha_sul1999, cha_sul2002,
coh2004}. This would follow, using the universality statement in
theorem A, from the existence of a theory of open-closed string
topology whose associated Calabi-Yau $A_\infty$ category was
equivalent to $H^{-\ast} (M)$.

Note that the degree shift in Chas-Sullivan's theory is incorporated
here in to the local system $\det^d$.

\subsection{The B model}
Let $X$ be a smooth projective Calabi-Yau variety of dimension $d$
over $\C$. Pick a holomorphic volume form on $X$.  Consider the dg
category $\op{Perf}(X)$, whose objects are bounded complexes of
holomorphic vector bundles on $X$, and whose morphisms are
$$
\op{Hom}_{\op{Perf}(X)}(E,F) = \Omega^{0,-\ast}( E^\vee \otimes F)
$$
(we reverse the grading, as all our differentials are homological).
The holomorphic volume form gives us a pairing
$$
\op{Hom}_{\op{Perf}(X)}(E,F) \otimes_\C \op{Hom}_{\op{Perf}(X)}(F,E)
\to \C
$$
of degree $d$, which is non-degenerate on homology.  Using the
homological perturbation lemma, we can transfer the $A_\infty$
structure on $\op{Perf}(X)$ to homology category. We should be able
to ensure that the resulting $A_\infty$ category is Calabi-Yau for
the natural pairing, using Hodge theory and the explicit form
of the homological perturbation lemma \cite{mer1999, kon_soi2001,  mar2004}. Denote by $\mc D^b_\infty(X)$ this Calabi-Yau $A_\infty$ category.

The closed TCFT $j^\ast \mbb{L} i_\ast \mc {D}^b_\infty (X)$ is the
$B$ model mirror to a TCFT constructed from Gromov-Witten invariants
of a compact symplectic manifold. We have seen that the homology of
$j^\ast \mbb{L} i_\ast  \mc {D}^b_\infty (X)$ is the Hochschild
homology of $ \mc {D}^b_\infty (X)$.

As the $A_\infty$ categories $\op{Perf}(X)$, $\mc D^b_\infty(X)$ are
quasi-isomorphic, they have the same Hochschild homology.  One
should be able to show that
\begin{equation}
HH_i(\mc D^b_\infty(X)) = HH_i(\op{Perf}(X)) =   \oplus_{q - p = i}
H^p(X,\Omega^q_X) \label{Hochschild Hodge}
\end{equation}
Theorem A, applied to $\mc D^b_\infty(X)$, implies there are
operations on $HH_\ast(\mc D^b_\infty(X))$ indexed by homology
classes on the moduli spaces of curves.  That is, if as before $\mc
M(I,J)$ is the moduli space of Riemann surfaces with $I$ incoming
and $J$ outgoing boundaries, there is a map
$$
H_\ast(\mc M(I,J),\det^d) \to \op{Hom} (HH_\ast(\mc
D^b_\infty(X))^{\otimes I}, HH_\ast(\mc D^b_\infty(X))^{\otimes J})
$$
compatible with gluing and disjoint union. These operations should
be the B-model mirror to corresponding operations on the homology of
a symplectic manifold coming from Gromov-Witten invariants.

Note that the usual derived category (without the $A_\infty$
enrichment) is a Calabi-Yau $A_\infty$ category. However, as usual,
passing to homology loses too much information.  This category
cannot encode the $B$ model.

\subsection{Gromov-Witten invariants and the Fukaya category}

The Fukaya category \cite{fuk_oh_oht_ono2000} of a symplectic
manifold should be an example of a unital Calabi-Yau $A_\infty$
category.  Thus, associated to the Fukaya category one has a closed
TCFT, whose homology is the Hochschild homology of the Fukaya
category.

Also, the Floer chains of the loop space of a symplectic manifold
should have a natural structure of closed TCFT, where the TCFT
operations come from counting pseudo-holomorphic maps. Thus to each
symplectic manifold we can associate two TCFTs, and it is natural to
conjecture that these are homotopy equivalent.  We will see that the
universality statement of theorem A allows us to relate these two
TCFTs, thus providing evidence for this conjecture.

\subsection{The TCFT associated to Gromov-Witten invariants}
First, let me explain a little about this second construction of a
TCFT, in the special case of a compact symplectic manifold $X$.  In
this case, the TCFT arises from Gromov-Witten invariants. Let
$\cmod$ be the Deligne-Mumford analog of Segal's category, that is
the category with objects finite sets, and morphisms stable
algebraic curves with incoming and outgoing marked points.  One can
find a homotopy equivalent model $\mc M'$ for Segal's category $\mc
M$ with a natural functor $\mc M' \to \cmod$. A chain-level theory
of Gromov-Witten invariants should give a functor from
$C_\ast(\cmod) \to \Comp$; pulling back via the functor $C_\ast(\mc
M') \to C_\ast(\cmod)$ will give the required TCFT. The model $\mc
M'$ for $\mc M$ we need was first constructed by Kimura, Stasheff
and Voronov in \cite{kim_sta_vor1995}. It can be constructed by
performing a real blow up of the Deligne-Mumford spaces along their
boundary. More precisely, we can take for $\mc M'$ the moduli space
of curves $\Sigma \in \cmod$, together with at each marked point a
section of the tautological $S^1$ bundle, and at each node a section
of the tensor product of the two tautological $S^1$ bundles
corresponding to either side of the node.

Suppose for simplicity that $c_1(X) = 0$, and let $\Sigma \in \mc
M(I,J)$. Then the real virtual dimension of the space of
pseudo-holomorphic maps from the fixed surface $\Sigma$ to $X$ is $d
( \chi(\Sigma) + \# I + \# J)$.  Thus, each such curve $\Sigma$
should give an operation
$$
C_\ast(X)^{\otimes I} \to C_\ast(X)^{\otimes J}
$$
of degree $d \chi + d \# J - d\# I$.  We want to construct a $d$
dimensional TCFT from a $2d$ dimensional symplectic manifold
\footnote{If $X$ does not satisfy $c_1(X) = 0$, we can work with
only a $\Z/2$ grading}. Therefore there should be a shift in
grading, and we should work with  $C_{\ast + d}(X)$.

One can check easily that if we work with this shift in grading, we
find a $d$ dimensional TCFT. The point is that the extra signs
arising from this shift in grading correspond to working with chains
on moduli space with coefficients in the local system $\op{det}^d$.

At the level of homology, this TCFT structure follows from the existence of
Gromov-Witten invariants; the chain level version we need is, I believe, still
conjectural.

\subsection{Comparing the TCFT associated to Gromov-Witten theory
with that from the Fukaya category}

Given a compact symplectic manifold, there should therefore be two
associated closed TCFTs: that coming from Gromov-Witten invariants,
and that constructed from the Fukaya category.  We now provide some
evidence for the conjecture that these are homotopy equivalent. Let
$X$ be a compact symplectic manifold of dimension $2d$, with Fukaya
category $\op{Fuk}(X)$.

\begin{conjecture}
There is a natural structure of $d$-dimensional open-closed TCFT,
whose D-branes are certain Lagrangian branes\footnote{Lagrangians
with the extra structure which makes them into an object of the
Fukaya category} in $X$, whose morphism spaces between D-branes
$L_1,L_2$ are the Lagrangian Floer chain groups
$$
\op{Hom}_i(L_1,L_2) = CF^{-i}( L_1,L_2)
$$
and whose complex of closed states is the shifted singular chain complex
$C_{\ast+d}(X)$ of $X$.
\end{conjecture}
This conjecture is I'm sure obvious to many people.  It is simply
asserting that the work of Fukaya-Oh-Ohta-Ono
\cite{fuk_oh_oht_ono2000} can be generalised to the case of Riemann
surfaces of all genus with open-closed boundary conditions, in a way
which takes into account families of surfaces.

Parts of this conjectural open-closed Gromov-Witten theory have
previously been constructed by P. Siedel \cite{sei2001,sei2001m} and
C.-C. Liu \cite{liu2002}. Seidel constructs the ``topological field
theory'' version with fixed complex structure on the source Riemann
surface. This corresponds to working with $H_0$ of moduli spaces.
The part dealing with only one Lagrangian, and varying source
Riemann surface, has been constructed by C.-C. Liu \cite{liu2002}.

A corollary of conjecture 1 and theorem  A is
\begin{ucorollary}
There is a map of closed TCFTs $j^\ast \mbb{L} i_\ast (\op{Fuk} X) \to
C_{\ast+d}(X)$ from the universal closed TCFT to the singular chains of $X$. On
homology this gives a map of homological TCFTs $HH_\ast(\op{Fuk} X) \to
H_{\ast+d}(X)$ from the Hochschild homology of the Fukaya category to the
homology of $X$.
\end{ucorollary}

A homological TCFT is like a TCFT except we replace the complex of chains on
moduli space by its homology.   The fact that the map $HH_\ast(\op{Fuk}(X)) \to
H_{\ast+d}(X)$ is a map of homological TCFTs means that it intertwines all
operations coming from the homology of moduli spaces of curves; that is the
diagram
$$
\xymatrix{ H_\ast(\mc M(I,J), \det^d) \otimes
HH_\ast(\op{Fuk}(X))^{\otimes I} \ar[r] \ar[d] &  H_\ast(\mc
M(I,J),\det^d) \otimes
H_{\ast+d}(X)^{\otimes I} \ar[d] \\
HH_\ast(\op{Fuk} (X))^{\otimes J} \ar[r] & H_{\ast+d}(X)^{\otimes J}}
$$
commutes.

The map from Hochschild to (Floer) homology is the same as that
constructed  by Seidel in \cite{sei2002}. The homology of a TCFT has
the structure of cocommutative coalgebra, coming from the
pair-of-pants coproduct. Note that as the pair of pants has Euler
characteristic $-1$, this is a map of degree $-d$. This coproduct
structure on $HH_\ast(\op{Fuk}(X))$ is dual to the standard cup
product on Hochschild cohomology, using the isomorphism
$HH_i(\op{Fuk}(X))^{\vee} \iso HH^{d + i }(\op{Fuk}(X))$.  The
coproduct on $H_{\ast+d}(X)$ is dual to the quantum cup product on
$H^\ast(X)$. Thus, the dual map $H^\ast(X) \to HH^\ast(\op{Fuk}(X))$
is in particular a ring homomorphism from quantum to Hochschild
cohomology.  Note that this dual map is of degree $0$.

Open-closed Gromov-Witten theory would give a map from the closed
TCFT associated to $\op{Fuk}(X)$ to that coming from the
Gromov-Witten theory of $X$.  It is natural to conjecture that this
is a quasi-isomorphism, that is
\begin{conjecture}
In good circumstances, the map $HH_\ast(\op{Fuk}(X)) \to
H_{\ast+d}(X)$ is an isomorphism.
\end{conjecture}
This conjecture, which was first proposed by Kontsevich
\cite{kon1995icm}, seems to be an integral part of the homological
mirror symmetry picture. Unfortunately, however, I really don't know
of much evidence.  Kontsevich presents a geometric motivation for
this conjecture in \cite{kon1995icm}, which I will reproduce here.
We can identify the Hochschild cohomology of the Fukaya category
with the endomorphisms of the identity functor, in the $A_\infty$
category of $A_\infty$ functors from $\op{Fuk}(X)$ to itself.  If we
could identify this $A_\infty$ category with $\op{Fuk}(X \times X,
\omega \oplus -\omega)$, as seems natural, we would see that the
Hochschild cohomology of $\op{Fuk}(X)$ would be the Lagrangian Floer
cohomology of the diagonal in $(X \times X,\omega \oplus -\omega)$,
which is known to coincide with the ordinary cohomology of $X$ with
the quantum product.

This conjecture implies that the homotopy Lie algebra controlling deformations
of $\op{Fuk}(X)$ is formal, and quasi-isomorphic to $H^\ast(X)$ with the
trivial Lie bracket. So that the formal neighbourhood of $\op{Fuk}(X)$ in the
moduli space of $A_\infty$ categories  is isomorphic to the formal
neighbourhood of the symplectic form in $H^\ast(X)$.  The homotopy Lie algebra
structure arises from an action of chains on moduli spaces of genus $0$ Riemann
surfaces. The homotopy Lie algebra structure on $C^\ast(X)$ should be trivial,
as the circle action is trivial.

\subsection{Acknowledgements}
I'd like to thank Ezra Getzler, Mike Hopkins, Andrey Lazarev, Tim Perutz, Graeme
Segal, Paul Seidel, Jim Stasheff, Dennis Sullivan, Constantin
Teleman, and Richard Thomas for  helpful conversations and
correspondence. I benefitted greatly from an inspiring talk of M.
Kontsevich at the Hodge centenary conference in 2003. In this talk
he described several results related to those in this paper; in
particular, he sketched a different construction of a TCFT structure
on the Hochschild chains of an $A_\infty$ algebra,  and also an
extension of this to the Deligne-Mumford spaces when the Hodge to de
Rham spectral sequence degenerates.

\subsection{Notation}
$\K$ will denote a field of characteristic zero. All homology and cohomology
will be with coefficients in $\K$, and all algebras and linear categories will
be defined over $\K$. $\Comp$ will denote the category of complexes of $\K$
vector spaces, with differential of degree $-1$, and with its standard
structure of symmetric monoidal category. For $r \in \Z$ we denote by $\K[r]$
the complex in degree $-r$, and for $V \in \Comp$ we write $V[r]$ for $V
\otimes \K[r]$. $\Vect$ will denote the category of $\Z$ graded $\K$ vector
spaces.

Instead of working with a field $\K$ and complexes of $\K$ vector spaces, the
main result remains true if instead we work with a commutative differential
graded algebra $R$ containing $\Q$, and flat dg $R$ modules.  (A dg $R$ module
$M$  is flat if the functor $M \otimes_R -$ is exact, that is takes
quasi-isomorphisms to quasi-isomorphisms).

\section{The open-closed moduli spaces in more detail}
\label{section oc moduli spaces}

Fix a set $\Lambda$ of D-branes.

A Riemann surface with open-closed boundary is a possibly disconnected Riemann
surface $\Sigma$, with boundary, some of whose boundary components are
parameterised in a way compatible with the orientation on $\Sigma$; these are
the incoming closed boundaries.   Other boundary components are parameterised
in the opposite sense; these are the outgoing closed boundaries.  There are
some disjoint intervals embedded in the remaining boundary components; these
are the open boundaries.  Some of these intervals are embedded in a way
compatible with the orientation on $\Sigma$; these are incoming open, the
remainder are outgoing open.

If we remove from $\partial \Sigma$ the open and closed
boundaries, what is left is a  one-manifold, whose connected
components are the free boundaries. Suppose the free boundaries of
$\Sigma$ are labelled by D-branes.  Then each open boundary $o$ of
$\Sigma$ has associated to it an ordered pair $(s(o),t(o))$ of
D-branes, associated to the free boundaries where it starts and
where it ends.

We require that each connected component of $\Sigma$ has at least
one incoming closed boundary or at least one free boundary.  We do
not impose a stability condition; note that no connected component
of $\Sigma$ can be a sphere or torus with no boundaries. However,
it is possible that a connected component of $\Sigma$ could be a
disc or an annulus with no open or closed boundaries, and only
free boundaries.  This would introduce an infinite automorphism
group; to remedy this, we replace the moduli space (stack) of
discs or annuli with no open or closed boundaries by a point.  One
can think of this as either taking the coarse moduli space, or
rigidifying in some way.

 Define a topological category $\mc M_\Lambda$. The objects of
$\mc M_\Lambda$ are quadruples $(C,O,s,t)$ where $C,O \in \Z_{\ge
0}$, and $s,t : O \to \Lambda$ are two maps.  (We use notation which
\linebreak  identifies the integer $O$ with the set
$\{0,1,\ldots,O-1\}$).  The space of morphisms
 $\mc M_\Lambda((C_+,O_+,s_+,t_+),(C_-,O_-,s_-,t_-))$ is the moduli space of
Riemann surfaces $\Sigma$ with open-closed boundary, with free
boundaries labelled by D-branes, with open incoming (respectively
outgoing) boundaries labelled by $O_+$ (respectively $O_-$), with
closed incoming (respectively outgoing) boundaries labelled by $C_+$
(respectively $C_-$), such that the maps $s_{\pm}, t_{\pm} : O_{\pm}
\to \Lambda$ coincide with those coming from the D-brane labelling
on $\Sigma$.  Composition in this category is given by gluing
incoming and outgoing open and closed boundaries to each other.

As defined, $\mc M_\Lambda$ is a non-unital category; it does not
have identity maps. To remedy this, we modify it a little. We
replace the moduli space of annuli, with one incoming and one
outgoing closed boundary, which is $\op{Diff}_+ S^1 \times_{S^1}
\op{Diff}_+ S^1 \times \R_{> 0}$, by the homotopy equivalent space
$\op{Diff}_+ S^1$, acting by reparameterisation. This should be
thought of as the moduli space of infinitely thin annuli. Similarly,
we replace the moduli space of discs with one incoming and one
outgoing open boundary by a point, which acts as the identity on the
open boundaries. We should perform this procedure also for any
surfaces which have connected components of one of these forms.

Disjoint union of surfaces and addition of integers $(C,O)$ makes
$\mc M_\Lambda$ into a symmetric monoidal topological category, in
the sense of \cite{mac1998}.  Note that this is a \emph{strict}
monoidal category; the monoidal structure is strictly associative.
It is not, however, strictly symmetric.

Let $C_\ast$ be the chain complex functor defined in the appendix,
from spaces to complexes of $\K$ vector spaces. $C_\ast$ is a symmetric monoidal functor, in the sense of \cite{mac1998}. This means that there is a
natural transformation $C_\ast(X) \otimes C_\ast(Y) \to C_\ast(X
\times Y)$, satisfying some coherence axioms. Define the category
$C_\ast(\mc M_\Lambda)$ to have the same objects as $\mc M_\Lambda$,
but with $C_\ast(\mc M_\Lambda)(l_1,l_2) = C_\ast(\mc
M_\Lambda(l_1,l_2))$ for $l_i \in \op{Ob} C_\ast(\mc M_\Lambda)$.
Since $C_\ast$ is a symmetric monoidal functor, $C_\ast(\mc M_\Lambda)$ is again
a symmetric monoidal category, but this time enriched over the
category of complexes. That is, $C_\ast(\mc M_\Lambda)$ is a
differential graded symmetric monoidal category.

As the set of D-branes will be fixed throughout the paper, we will
occasionally omit the subscript $\Lambda$ from the notation.
\begin{definition}
Let $\mscr{OC}_\Lambda = C_\ast(\mc M_\Lambda)$.  Let $\mscr{O}_\Lambda$ be the
full  subcategory whose objects are $(0, O, s,t)$, that is have no closed part.
Let $\mscr{C}$ be the subcategory whose objects have no open part, and whose
morphisms are Riemann surfaces with only closed boundaries.   $\mscr{C}$ is
independent of $\Lambda$. These categories are differential graded symmetric
monoidal categories.
\end{definition}

Note that if $\Lambda \to \Lambda'$ is a map of sets, there are
corresponding functors $\mscr{OC}_\Lambda \to
\mscr{OC}_{\Lambda'}$ and $\mscr{O}_\Lambda \to
\mscr{O}_{\Lambda'}$. We could think of $\mscr{O}$ and $\mscr{OC}$
as categories fibred over the category of sets.

We need twisted versions of these categories. Consider the graded
$\K$ local system $\det$ on the spaces of morphisms in $\mc
M_\Lambda$, whose fibre at a surface $\Sigma$ is
\begin{align*}
\det(\Sigma) &=  \op{det} ( H^0(\Sigma) - H^1(\Sigma) + \K^{O_-} ) [ O_- -
\chi(\Sigma) ] \\ &= \op{det}( H^0(\Sigma) \oplus \K^{O_-} ) \otimes (\op{det}
H^1(\Sigma))^{-1} [O_- - \chi(\Sigma)]
\end{align*}
Here $O_-$ is the number of open outgoing boundary components of $\Sigma$. The
number in square brackets refers to a shift of degree; so this is a graded
local system situated in degree $\chi(\Sigma) - O_-$.

Suppose $\Sigma_1, \Sigma_2$ are composable morphisms in $\mc
M_\Lambda$. Then there is a natural isomorphism
$$
\op{det} \Sigma_2 \otimes \op{det} \Sigma_1 \to \op{det}( \Sigma_2 \circ
\Sigma_1)
 $$
This follows from the Mayer-Vietoris exact sequence obtained from
writing $\Sigma_2 \circ \Sigma_1$ as a union of the $\Sigma_i$'s.
Let $C_-^i$ and $O_-^i$ be the open and closed outgoing boundaries
of $\Sigma^i$. We have
\begin{multline*}
0 \to H^0 (\Sigma_2 \circ \Sigma_1) \to H^0(\Sigma_2) \oplus H^0(\Sigma_1) \to
\K^{C_-^1 + O_-^1} \\ \to H^1(\Sigma_2 \circ \Sigma_1) \to H^1(\Sigma_2) \oplus
H^1(\Sigma_1) \to \K^{C_-^1} \to 0
\end{multline*}
Here $\K^{C_-^1 + O_-^1}$ arises as $H^0( \Sigma_2 \cap \Sigma_1)$ and
$\K^{C_-^1}$ arises as $H^1(\Sigma_2 \cap \Sigma_1)$.  Note that the
orientation on the outgoing boundary of $\Sigma_1$ gives a natural isomorphism
$H^1(\Sigma_2 \cap \Sigma_1)\iso \K^{C_-^1}$.

We will see  that $(\mc M_\Lambda,\det)$ again forms a kind of
category. Consider the symmetric monoidal category whose objects are
pairs $(X,E)$ where $X$ is a  topological space and $E$ is a graded
$\K$ local system on $X$, such that a map $(X,E) \to (Y,F)$ is a map
$f : X \to Y$ and a map $E \to f^\ast F$, and such that
$$(X,E) \otimes (Y,F) = (X \times Y, \pi_1^\ast E \otimes \pi_1^\ast F)$$
The symmetrisation isomorphism $(X,E)\otimes(Y,F) \iso (Y,F)
\otimes (X,E)$ as usual picks up signs from the grading on $E$ and
$F$.

We want to show that $(\mc M_\Lambda,\det)$ forms a category
enriched over the category of spaces with graded local systems.  All
that needs to be checked is that for composable surfaces
$\Sigma_1,\Sigma_2,\Sigma_3$ the diagram
$$
\xymatrix{ \op{det} (\Sigma_3) \otimes \op{det} (\Sigma_2) \otimes \op{det}
(\Sigma_1) \ar[r] \ar[d] & \op{det} (\Sigma_3) \otimes
\op{det} (\Sigma_2 \circ \Sigma_1) \ar[d] \\
\op{det} (\Sigma_3 \circ \Sigma_2) \otimes \op{det} (\Sigma_1) \ar[r] &
\op{det} (\Sigma_3 \circ \Sigma_2 \circ \Sigma_1) }
$$
commutes, where $\op{det} (\Sigma)$ is the fibre of the local
system at $\Sigma$.  This is a fairly straightforward calculation.

There is also a natural isomorphism $\op{det}( \Sigma_1 \amalg
\Sigma_2) \iso \op{det} (\Sigma_1) \otimes \op{det} (\Sigma_2)$.
This gives $(\mc M_\Lambda,\det)$ the structure of symmetric
monoidal category.

The functor $C_\ast$ defined in the appendix is a functor from the
category of spaces with graded $\K$ local systems to complexes,
which computes homology with local coefficients.  Since $C_\ast$ is
a symmetric monoidal functor, it follows that $C_\ast(\mc M_\Lambda,\det)$ is
again a symmetric monoidal category.

We can think of the chain category $C_\ast(\mc M_{\Lambda},\det)$
geometrically as follows.  A chain with local coefficients on
$(X,E)$ can be thought of as a singular  simplex $f : \diag_n \to X$
together with a section of $f^\ast E \otimes \omega$, where $\omega$
is the orientation sheaf on $\diag_n$. Thus a chain in $C_\ast(\mc
M_{\Lambda},\det)$ should be thought of as an oriented $n$ parameter
family of Riemann surfaces $\Sigma$ with a section of
$\det(\Sigma)$.

We can also twist $\mc M_{\Lambda}$ by tensor powers $\det^d =
\det^{\otimes d}$, where $d \in \Z$.
\begin{definition}
Let $\mscr{OC}_\Lambda^d$ be the category $C_\ast(\mc M_\Lambda,\det^d)$.  As
before, let $\mscr{O}_\Lambda^d$ be the full subcategory whose objects have no
closed part, and let $\mscr{C}^d$ be the  subcategory whose objects have no
open part and whose morphisms have only closed boundaries.  These are
differential graded symmetric monoidal categories.
\end{definition}
As before, if $\Lambda \to \Lambda'$ is a map of sets, there are
corresponding functors $\mscr{O}^d_\Lambda \to
\mscr{O}^d_{\Lambda'}$ and $\mscr{OC}^d_\Lambda \to
\mscr{OC}^d_{\Lambda'}$.

\section{Some homological algebra for symmetric monoidal categories}
\label{general stuff}

\subsection{Differential graded symmetric monoidal categories}

We work with differential graded symmetric monoidal categories,
over $\K$. Symmetric monoidal is in the sense of  MacLane
\cite{mac1998}; differential graded means that the morphism spaces
are complexes of $\K$ vector spaces (with differential of degree
$-1$), and the composition maps are bilinear and compatible with
the differential.  Call these dgsm categories, for short. A good
reference for the general theory of dg categories is
\cite{kel1994}.

The dgsm categories controlling topological conformal field theory
are strictly monoidal.  On objects, $(\alpha \amalg \beta) \amalg
\gamma  = \alpha \amalg (\beta \amalg \gamma)$, and similarly the
diagram
$$
\xymatrix{ \op{Hom} ( \alpha,\alpha') \otimes \op{Hom}
(\beta,\beta') \otimes \op{Hom} ( \gamma,\gamma' ) \ar[r] \ar[d] &
\op{Hom} ( \alpha \amalg \beta,\alpha' \amalg \beta') \otimes
\op{Hom} ( \gamma , \gamma' ) \ar[d] \\
\op{Hom} ( \alpha ,\alpha') \otimes \op{Hom}( \beta \amalg \gamma,
\beta' \amalg \gamma' ) \ar[r] &  \op{Hom}( \alpha \amalg \beta
\amalg \gamma, \alpha' \amalg \beta' \amalg \gamma' ) }
$$
commutes. (We use $\amalg$ and $\otimes$ interchangeably for the tensor product
in the categories controlling TCFT). However, the symmetry isomorphism $\alpha
\amalg \beta \to \beta \amalg \alpha$ is not an identity, nor do we always have
$\alpha \amalg \beta = \beta \amalg \alpha$. If $A$ is strictly monoidal, for
each $\sigma \in S_n$ there is an isomorphism $a_1 \otimes \ldots a_n \iso
a_{\sigma(1)} \otimes \ldots \otimes a_{\sigma(n)}$, compatible with
composition in the symmetric groups.

Let $A,B$ be dgsm categories, which for simplicity we assume are
strictly monoidal.  A monoidal functor $F : A \to B$ is a functor
$F$, compatible with the dg structure, together with natural
transformations $F(a)  \otimes F(a') \to F(a \otimes a')$, such
that the diagrams
$$
\xymatrix{ F (a) \otimes F(a') \otimes F(a'') \ar[r] \ar[d] & F(a
\otimes a') \otimes F(a'') \ar[d]  \\
F(a) \otimes F(a' \otimes a'') \ar[r] & F(a \otimes a' \otimes
a'') }
$$
and
$$
\xymatrix{ F(a) \otimes F(a') \ar[r] \ar[d] &  F(a \otimes a')
\ar[d] \\
F(a') \otimes F(a) \ar[r] &   F(a' \otimes a) }
$$
commute.

Although our dgsm categories may (or may not) have an object which
is a unit for the tensor product, we do not assume the functor $F$
takes units to units.

To a dgsm category $A$ are associated several important auxiliary
categories. First there is the homology category $H_\ast A$, whose
objects are the same as those of $A$, but with
$$
\op{Hom}_{H_\ast A}(a,a') = H_\ast \op{Hom}_A(a,a')
$$
$H_\ast A$ is a graded symmetric monoidal category; the morphisms
are graded vector spaces.  Similarly we have the category $H_0 A$,
whose morphisms are $H_0 \op{Hom}_A(a,a')$.  Also, there is the
category $Z_0 A$, which is a subcategory of $A$ with the same
objects, but whose morphisms are closed maps  of degree $0$. A map
in $Z_0 A$ is called a quasi-isomorphism if it is an isomorphism
in $H_0 A$.

One example of a dgsm category is the category $\Comp$ of
complexes of $\K$ vector spaces. The monoidal structure is given
by tensor product.

A left $A$ module is a (monoidal) functor $A \to \Comp$.  A right
$A$ module is a (monoidal) functor $A^{op} \to \Comp$, where
$A^{op}$ is the opposite category to $A$. If $M,N : A \to B$ are
monoidal functors to a dgsm category $B$, a natural transformation
$\phi : M \to N$ consists of a collection of maps $\phi(a) \in
\op{Hom}_B( M(a), N(a))$ satisfying the following conditions.
\begin{enumerate}
\item
$\phi(a)$ is natural for morphisms in $a$.  That is, if $f : a \to
a'$ then $\phi(a') M(f) = N(f) \phi(a)$.
\item
The morphisms $\phi(a) \in \op{Hom}_B(M(a),N(a))$ are all closed
and of degree $0$.
\item
The diagram
$$
\xymatrix{M(a) \otimes M(a') \ar[r] \ar[d] & N(a) \otimes N(a')
\ar[d] \\
M(a \otimes a') \ar[r] & N(a \otimes a')}
$$
commutes.
\end{enumerate}
Thus for example we have a category $A-\mod$ of left $A$ modules
and $\mod-A$ of right $A$ modules.  Note that $A-\mod$ is just a
category, not a dg category; it is not even an additive category.

If $A,B$ are dgsm categories, we can form their tensor product $A
\otimes B$.  The objects are
$$
\op{Ob} (A \otimes B) = \op{Ob} A \times \op{Ob} B
$$
and the morphisms are described by
$$
\op{Hom}( a \times b, a' \times b') = \op{Hom}(a,a') \otimes_\K
\op{Hom}(b,b')
$$
$A \otimes B$ is again a dgsm category.  An $A-B$ bimodule is a
monoidal functor $A \otimes B^{op} \to \Comp$.

We will often use the notation
$$
A(a,a') = \op{Hom}_A(a,a')
$$
$A$ defines an $A-A$ bimodule over itself, by the functor $A \otimes A^{op} \to
\Comp$ which sends
$$
(a_1,a_2) \mapsto A(a_2,a_1)
$$
However, if $a \in \op{Ob} A$ is an object, then the functor $A \to \Comp$
defined by $\op{Hom}(a,-)$ is not in general monoidal, and so does not give an
$A$-module in our sense.

\subsection{Notation about exact functors}

Suppose a category $C$ has a notion of quasi-isomorphism.  That
is, suppose we are given a subset of the set of morphisms of $C$,
which is closed under composition and which contains all
isomorphisms. We say objects in $C$ are quasi-isomorphic if they
can be connected by a chain of quasi-isomorphisms. We write $c
\simeq c'$ to indicate that $c,c'$ are quasi-isomorphic.

If $D$ also has a class of quasi-isomorphisms, a functor $F : C \to D$ is
called exact if it takes quasi-isomorphisms to quasi-isomorphisms.

A natural transformation between exact functors $F,G : C \to D$ is called a
quasi-isomorphism if for each $c \in C$ the maps $F(c) \to G(c)$ are
quasi-isomorphisms.
\begin{definition}
A quasi equivalence between $C,D$ is a pair of functors $F : C \to D$ and $G :
D \to C$ such that $F \circ G$ is quasi-isomorphic to $\op{Id}_D$, and $G \circ
F$ is quasi-isomorphic to $\op{Id}_C$.  That is,
\begin{align*}
F \circ G &\simeq \op{Id}_D & G \circ F & \simeq \op{Id}_C
\end{align*}
\label{def quasi_equiv}
\end{definition}

For example, let $A$ be a dgsm category. Recall $Z_0 A$ is the
category with the same objects as $A$ but whose morphisms are
closed of degree $0$.  A map $a \to a'$ in $Z_0 A$ is a
quasi-isomorphism if it is an isomorphism in $H_0 A$.

Any functor $F : A \to B$ between dgsm categories restricts to an
exact functor $Z_0 A \to Z_0 B$.  Thus the category of functors $A
\to B$ acquires a notion of quasi-isomorphism. In particular we
can talk about quasi-isomorphisms in $A-\mod$; these are just
morphisms which are quasi-isomorphisms of the underlying
complexes.

We would like to do some kind of homotopy theory with categories
$A-\mod$ for various $A$.  I am going to do this in  a slightly ad
hoc fashion.  Probably one should put some extra structure on the
categories $A-\mod$ which would allow a more canonical notion of
derived functor.  For example, one could try to make $A-\mod$ into
a closed model category whose weak equivalences are
quasi-isomorphisms.  However, closed model structures are
difficult to construct.  One alternative structure which seems
more natural in this situation would be to consider categories
fibred over the category of differential graded commutative
algebras.  The fibre over $R$ should be the category of $R$ linear
functors from $A \otimes R$ to flat complexes of $R$ modules.  One
could use this structure to define notions of homotopy between
maps, and eventually to define derived functors in a more
canonical way.

Instead of attempting to construct any such general theory, or
give a closed model structure on $A-\mod$,  I will perform
homotopic constructions in a slightly more ad hoc fashion.  We
only need to derive one kind of functor; if $f : A \to B$ is a
functor, there is a pull back functor $f^\ast : B-\mod \to
A-\mod$, which is exact, and a left adjoint $f_\ast$ which is not.
We will construct the derived functor of $f_\ast$;  it will be
clear from the construction that there is a unique left derived
functor $\mbb{L}f_\ast$ up to quasi-isomorphism.

\subsection{Derived tensor products}
Let $M$ be a $B-A$ bimodule.  Let $N$ be a left $A$ module.  Then
we can form a left $B$ module $M \otimes_A N$. For each $b \in B$,
$M \otimes_A N(b)$ is defined to be the universal complex with
maps $M(b,a) \otimes_\K N(a) \to (M \otimes_A N)(b)$, such that
the diagram
$$
\xymatrix{  M(b,a) \otimes_\K A(a',a) \otimes_\K N(a') \ar[r]
\ar[d] & M(b,a) \otimes_\K N(a) \ar[d] \\
M(b,a') \otimes_\K N(a') \ar[r] & (M \otimes_A N)(b) }
$$
commutes.  One can check that $M \otimes_A N$ is again a monoidal
functor from $B$ to complexes.  Thus $M \otimes_A -$ defines a
functor $A-\mod \to B-\mod$.

Let $f : A \to B$ be a functor between dgsm categories.  Then $B$
is a $B-B$ bimodule, and so becomes an $A-B$ bimodule and a $B-A$
bimodule via the functors $A\otimes B^{op} \to B \otimes B^{op}$
and $B \otimes A^{op} \to B \otimes B^{op}$.  We can define
functors $f_\ast : A-\mod \to B-\mod$ by
$$f_\ast M = B \otimes_A M$$
and $f^\ast : B-\mod \to A-\mod$ by defining $f^\ast N$ to be $N$ with the
induced $A$ action.  So as a functor $A \to \Comp$,  $f^\ast N$ is the
composition of $N : B \to \Comp$ with the functor $f : A \to B$.   The push
forward functor $f_\ast$ is the left adjoint to $f^\ast$.

Note that $f^\ast$ is exact.  In general $f_\ast$ is not exact.
However, we can construct a derived version of $f_\ast$ which is
exact.

We say an  $A$ module $M$ is flat if the functor $- \otimes_A M$ from right $A$
modules to complexes is exact. Let $A-\op{flat}$ be the subcategory of flat $A$
modules, and let $i : A-\op{flat} \into A -\mod$ be the inclusion.

\begin{definition}
Let $\op{SymOb} A \subset A$ be the subcategory with the same
objects, but whose morphisms are the identity maps and the symmetry
isomorphisms $a_1 \otimes \ldots a_n \iso a_{\sigma(1)} \otimes
\ldots \otimes a_{\sigma(n)}$, for $\sigma \in S_n$. $\op{SymOb} A$
is again a symmetric monoidal category with a monoidal functor
$\op{SymOb} A \to A$.   Also $\op{SymOb} A$ is a groupoid.

Let $\op{SymOb}_\K A \subset A$ be the sub linear category whose
morphisms are spanned by those of $\op{SymOb} A$.
\end{definition}

\begin{theorem}
Suppose $A$ is strictly monoidal, $A$ has a unit for the tensor
product, and suppose the groupoid $\op{SymOb} A$ has finite
automorphism groups for any object. Then there is a functor $F :
A-\mod \to A-\op{flat}$ such that $F \circ i$ and $i \circ F$ are
quasi-isomorphic to the identity functors. \label{theorem flat
resolution}
\end{theorem}
The conditions of the theorem hold for the dgsm categories
controlling topological conformal field theory.

This result is false except in characteristic zero. We will assume
the conditions of the theorem for the dgsm categories $A$, $B$ we
use for the rest of this section.

An $A-B$ bimodule $M$ is called $A$--flat if the functor $ -
\otimes_A M$ is exact, as a functor from right $A$ modules to
right $B$ modules.  The proof of this result will also show that
there exists functorial $A$--flat resolutions of $A-B$ bimodules.

\begin{proof}
Let $\op{Ord}$ be the simplicial category, whose objects are the
non-empty totally ordered finite sets, and whose morphisms are
non-decreasing maps. We will refer to the object $\{1,\ldots,n\}$ of
$\op{Ord}$ by $\{n\}$.  For a category $\mc C$, a simplicial object
of $\mc C$ is a functor $\op{Ord}^{\op{op}} \to \mc C$.  If $M$ is a
simplicial object of $\mc C$, we will  write $M\{n\} \in \op{Ob} \mc
C$ for the $n$ simplices of $M$.

For each $A$ module $M$, define a simplicial $A$ module
$\op{Bar}_A^\triangle M$ to have for $n$ simplices the $A$ module
$$
\op{Bar}_A^\triangle M \{n\} =  A \otimes_{\op{SymOb}_\K A} A
\otimes \ldots \otimes_{\op{SymOb}_\K A} A \otimes M =
A^{\otimes_{\op{SymOb}_\K A} n}\otimes_{\op{SymOb}_\K A} M
$$
The face maps come from the product maps of $\op{SymOb}_\K A$
bimodules $A \otimes_{\op{SymOb}_\K A} A \to A$, and the map $A
\otimes_{\op{SymOb}_\K} A M \to M$ of left $\op{SymOb}_\K A$
modules. The degeneracy maps come from the map $\op{SymOb}_\K A \to
A$ of $\op{SymOb}_\K A$ bimodules.

Denote by $\op{Comp}_\K^\tr$ the category of simplicial chain
complexes, that is functors $\op{Ord}^{\op{op}} \to \Comp$.  This is
a symmetric monoidal category.  The tensor product is pointwise; if
$C,D$ are simplicial chain complexes, then $(C \otimes D)\{n\} =
C\{n\} \otimes D[n]$.   A simplicial $A$ module is the same as a
symmetric monoidal functor $A \to \Comp^\tr$.

The normalised realisation functor $\abs{\quad} : \Comp^\tr \to
\Comp$ is defined by sending a simplicial chain complex $C$ to
$$\abs{C} = \oplus_{n >
0} C\{n\} / C^{degenerate}\{n\}[-n]$$ Here $C^{degenerate}\{n\}$ is
the image of the degeneracy maps. The symbol $[-n]$ refers to a
shift in degree. The differential on $\abs{C}$ is composed of the
differential on the summands $C\{n\}/ C^{degenerate}\{n\}$ and the
alternating sum of the face maps.

A map $C \to D$ of simplicial chain complexes is a quasi-isomorphism
if the maps $C\{n\} \to D\{n\}$ are quasi-isomorphisms.  The
realisation functor $\abs{\quad}$ is exact, that is it takes
quasi-isomorphisms to quasi-isomorphisms.

The shuffle product maps
$$
\abs{C} \otimes \abs{D} \to \abs{C \otimes D}
$$
make $\abs{\quad}$ into a symmetric monoidal functor.

Thus, in particular,
$$
\op{Bar}_A M \defeq \abs{\op{Bar}_A^\triangle M }
$$
is a symmetric monoidal functor $A \to \Comp$, in other words a left
$A$ module.

We can consider $M$ as a constant simplicial $A$ module.  There is a
natural map $\op{Bar}_A^\tr M \to M$, which on $n$ simplices comes
from the product map $A^{\otimes_{\op{SymOb}_\K} n}
\otimes_{\op{SymOb}_\K A} M \to M$. This induces a map of
realisations $\op{Bar}_A M \to M$.
\begin{lemma}
The map $\op{Bar}_A M \to M $ is a quasi-isomorphism.
\end{lemma}
\begin{proof}
This is essentially standard.
\end{proof}

\begin{lemma}
For any $A$ module $M$, $\op{Bar}_A M$ is a flat $A$ module.
\end{lemma}
\begin{proof}
Let $N,N'$ be right $A$ modules, with a quasi-isomorphism $N \to
N'$.  We need to show that the map
$$
N \otimes \op{Bar}_A M \to N' \otimes \op{Bar}_A M
$$
is a quasi-isomorphism.

We can consider $N,N'$ as constant simplicial right $A$ modules, and
form the tensor product simplicial chain complex $N \otimes_A
\op{Bar}^\tr_A M$.  This has for $n$ simplices
$$
N \otimes_A \op{Bar}^\tr_A M \{n\} = N \otimes_{\op{SymOb}_\K A}
A^{\otimes_{\op{SymOb}_\K A} n-1} \otimes_{\op{SymOb}_\K A} M
$$
It is easy to see that $N \otimes \op{Bar}_A M$ is the realisation
of this simplicial chain complex.  Since the realisation functor is
exact, it suffices to show that the map
$$
N \otimes \op{Bar}^\tr_A M \to N' \otimes \op{Bar}^\tr_A M
$$
of simplicial chain complexes is a quasi-isomorphism. To show this,
it suffices to show that the map
$$
N \otimes_{\op{SymOb}_\K A} A^{\otimes_{\op{SymOb}_\K A} n-1}
\otimes_{\op{SymOb}_\K A} M
 \to   N'
\otimes_{\op{SymOb}_\K A} A^{\otimes_{\op{SymOb}_\K A} n-1}
\otimes_{\op{SymOb}_\K A} M
$$
is a quasi-isomorphism.  More generally, if $P$ is any left
$\op{SymOb}_\K A$ module, consider the map
$$
N \otimes_{\op{SymOb}_\K A} P \to N' \otimes_{\op{SymOb}_\K A} P
$$
This is always a quasi-isomorphism.  As, tensor product over
$\op{SymOb}_\K A$ only involves taking coinvariants for finite group
actions (using the assumption that all isomorphism groups in the
groupoid $\op{SymOb} A$ are finite).  As we are working in
characteristic zero, the functor of coinvariants for a finite group
action is exact.

\end{proof}
This completes the proof of  theorem \ref{theorem flat resolution}.
A similar argument shows that there exists functorial $A$--flat
resolutions for $A-B$ bimodules.
\end{proof}

If $M$ is an $A-B$ bimodule, and $N$ is a left $B$ module, define a left $A$
module by
$$M \otimes_B^{\mbb L} N = M \otimes_B \op{Bar}_B N$$
Any other flat resolution of $N$ will give a quasi-isomorphic
answer; as, suppose $N',N''$ are flat resolutions of $N$, and $M'$
is a $B$ flat resolution of $M$.  Then
$$M \otimes_B N' \simeq M' \otimes_B N'
\simeq M' \otimes_B N'' \simeq M \otimes_B N''$$
\subsection{Derived push forwards}
Let $f : A \to B$ be a functor, between dg symmetric monoidal
categories. Let $N$ be a left $A$ module. Define
$$
\mbb {L} f_\ast N =  B \otimes_A^{\mbb L} N
$$
Note that if we define $\mbb{L}' f_\ast N = B ' \otimes_A N$,
where $B'$ is an $A$--flat resolution of $B$, then $\mbb{L}'
f_\ast$ and $\mbb{L} f_\ast$ are quasi-isomorphic functors. Also,
if we took any other choice of functorial flat resolution of $N$
we would get a quasi-isomorphic left derived functor.

Recall that $f^\ast$ is defined by considering a left $B$ module
as a left $A$ module.
\begin{lemma}
Suppose that the homology functor $H_\ast(f) : H_\ast(A) \to
H_\ast(B)$ is fully faithful.  Then the functor $f^\ast \mbb{L}
f_\ast$ is quasi-isomorphic to the identity functor on $A-\mod$.
\end{lemma}
\begin{proof}
Note that
$$
f^\ast \mbb{L} f_\ast N = B \otimes_{A} \op{Bar}_A N
$$
considered as a left $A$ module.  There is a map $A \to B$ of $A$
bimodules, and so a map
$$\op{Bar}_A N = A \otimes_{A} \op{Bar}_A N \to B \otimes_{A} \op{Bar}_A
N$$ We need to show this is a quasi-isomorphism.  This is clear,
as $\op{Bar}_A N$ is flat as an $A$ module, and the the statement
that $H_\ast(f)$ is full and faithful means that the map $A \to B$
of $A-A$ bimodules is a quasi-isomorphism.
\end{proof}

\begin{definition}
A quasi-isomorphism between dgsm categories is a functor $f : A
\to B$ such that $H_\ast(f)$ is full and faithful and $f$ induces
an isomorphism on the set of objects.
\end{definition}

\begin{theorem}
If $f : A \to B$ is a quasi-isomorphism, then the functors
$\mbb{L} f_\ast$ and $f^\ast$ are inverse quasi-equivalences
between $A-\mod$ and $B-\mod$, and between $\mod-A$ and $\mod-B$.
\label{theorem quasi equivalence module}
\end{theorem}
\begin{proof}
 We have seen that the functor $f^\ast
\mbb{L} f_\ast$ is quasi-isomorphic to the identity functor on $A-\mod$.  We
need to show that $\mbb{L} f_\ast f^\ast$ is quasi-isomorphic to the identity
functor on $B-\mod$.

Note that
$$
\mbb{L} f_\ast f^\ast N = B \otimes_A^{\mbb L} N = B \otimes^{\mbb
L}_A B \otimes_B N
$$
Therefore it suffices to write down a weak equivalence $B
\otimes_A^{\mbb L} B \to B$ of $B-B$ bialgebras.

The $B-B$ bialgebra $B \otimes_A^{\mbb L} B$ is the realisation of
the simplicial $B-B$ bialgebra $B \otimes_A \op{Bar}^\triangle_A B$,
using the notation of the proof of theorem \ref{theorem flat
resolution}.

By assumption, the functor $A \to B$ induces an isomorphism on the
set of objects.  Thus the categories $\op{SymOb}_\K A$ and
$\op{SymOb}_\K B$ are isomorphic; let us use the notation
$$
C \defeq \op{SymOb}_\K A = \op{SymOb}_\K B
$$
The $n$ simplices of $B \otimes_A \op{Bar}^\triangle_A B$
$$
B \otimes_{C} A^{\otimes_C n-1} \otimes_C B
$$
The map $A \to B$ of $C$ bimodules induces a map
$$
B \otimes_{C} A^{\otimes_C n-1} \otimes_C B \to B \otimes_{C}
B^{\otimes_C n-1} \otimes_C B
$$
which is a quasi-isomorphism, because tensor product of $C$
bimodules is an exact functor.  It is easy to see that this map is
the $n$th component of a quasi-isomorphism of simplicial chain
complexes
$$B \otimes_A \op{Bar}_A^\tr B \to B \otimes_B \op{Bar}_B^\tr B$$
The associated map on realisations is a quasi-isomorphism.  There is
a natural quasi-isomorphism of $B-B$ bimodules $B \otimes_{B}
\op{Bar}_B B \to B$.  Putting these quasi-isomorphisms together we
get a quasi-isomorphism $B \otimes^{\mbb L}_A B \to B$.

\end{proof}

\begin{lemma}
\label{lemma tensor compatible} Denote also by $f^\ast$ and $\mbb
L f_\ast$ the induced quasi-equivalences $\mod - A \times A- \mod
\leftrightarrows \mod-B \times B -\mod$. Both triangles in the
following diagram commute up to quasi-isomorphism.
$$
\xymatrix{ \mod - A \times A -\mod \ar@<1ex>[dd] ^{\mbb{L} f_\ast}
\ar
[drr]^{\otimes^{\mbb L}} \\
& & \Comp \\
 \mod - B \times B - \mod \ar@<1ex>[uu]^{f^\ast} \ar[urr]
^{\otimes^{\mbb L}} }
$$
The diagonal arrows are the tensor product maps which take a pair
$(M,N)$ where $M \in \mod-A$ and $N \in A-\mod$ to $M
\otimes^{\mbb L}_A N$.
\end{lemma}
\begin{proof}
It is sufficient to show that one of the triangles commutes up to
quasi-isomorphism. So we need to show that
$$
\mbb{L} f_\ast (M) \otimes^{\mbb L} _ B \mbb{L} f_\ast (N) \simeq
M \otimes_A^{\mbb L} N
$$
This follows from the chain of quasi-isomorphisms
\begin{align*}
\mbb{L}f_\ast M \otimes^{\mbb L}_B \mbb {L} f_\ast N & \simeq M
\otimes_A^{\mbb L} B \otimes_B^{\mbb L} B \otimes_A^{\mbb L} N \\
&\simeq M \otimes_A^{\mbb L} ( B \otimes_{A}^{\mbb L} N) \\
& \simeq M \otimes_A^{\mbb L} N
\end{align*}
where the last quasi-isomorphism comes from the quasi-isomorphism $N \simeq
f^\ast \mbb{L} f_\ast N = B \otimes_{A}^{\mbb L} N$ as an $A$-module.
\end{proof}
\section{Topological conformal field theories}

A symmetric monoidal functor $F: A \to B$ between dgsm categories is called
\emph{split} if the maps $F(a) \otimes F(a') \to F(a\otimes a')$
are all isomorphisms. This is what MacLane \cite{mac1998} calls
strong.   $F$ is called h-split, or homologically split, if
$H_\ast(F) : H_\ast(A) \to H_\ast(B)$ is split.  Note that being
h-split is an exact condition : if $F \simeq F'$ then $F$ is
h-split if and only if $F'$ is.
\begin{definition}
\begin{enumerate}
\item
An open topological conformal field theory  of dimension $d$ is a pair
$(\Lambda, \Phi)$ where $\Lambda$ is a set of D-branes, and $\Phi \in
\mscr{O}^d_{\Lambda}-\mod$ is a symmetric monoidal functor
$$\Phi : \mscr{O}^d_{\Lambda} \to \Comp$$
which is h-split.

A morphism of open TCFTs $(\Lambda,\Phi) \to (\Lambda',\Phi')$  is
a map $\Lambda \to \Lambda'$ of sets, and a morphism $\Phi \to
f^\ast \Phi'$ in $\mscr{O}^d_{\Lambda}-\mod$. Here $f :
\mscr{O}^d_{\Lambda} \to \mscr{O}^d_{\Lambda'}$ is the functor
induced by the map of sets $\Lambda \to \Lambda'$.
\item
A closed topological conformal field theory of dimension $d$ is a h-split
symmetric monoidal functor $\mscr{C}^d \to \Comp$.  A morphism of closed TCFTs is a
morphism in $\mscr{C}^d -\mod$.
\item
An open-closed topological conformal field theory of dimension $d$ is a pair
$(\Lambda,\Phi)$ where $\Lambda$ is a set of D-branes and $\Phi$ is a symmetric monoidal functor
$$\Phi : \mscr{OC}^d_{\Lambda} \to \Comp$$
which is h-split.

A morphism of open-closed TCFTs $(\Lambda,\Phi) \to
(\Lambda',\Phi')$ is a map $\Lambda \to \Lambda'$ of sets, and a
morphism $\Phi \to f^\ast \Phi'$ in $\mscr{OC}^d_{\Lambda}-\mod$.
Here $f : \mscr{OC}^d_{\Lambda} \to \mscr{OC}^d_{\Lambda'}$ is the
functor induced by the map of sets $\Lambda \to \Lambda'$.
\end{enumerate}
\end{definition}
The condition that the functors are h-split is important.  For
example, if $\Psi$ is a closed TCFT, then this means that
$$
H_\ast(\Psi(C)) = H_\ast(\Psi(1))^{\otimes C}$$ where $C$ is the
number of closed boundaries. Thus, if $\Psi$ is a a closed TCFT we
can talk about its homology, which is just a graded vector space;
we will use the notation $H_\ast(\Psi)$ for $H_\ast(\Psi(1))$.
Then $H_\ast(\Psi)$ carries operations from the homology of moduli
spaces of curves.  That is, there are maps
$$
H_\ast(\mscr{C}^d(I,J))  \to \op{Hom} ( H_\ast(\Psi)^{\otimes I},
H_\ast(\Psi)^{\otimes J})
$$

A pair $\lambda_1,\lambda_2$ of D-branes gives an object
$\{\lambda_1,\lambda_2\}$ of $\mscr{O}_\Lambda^d$, corresponding to one open
boundary from $\lambda_1$ to $\lambda_2$. For an open TCFT $(\Lambda,\Phi)$ we
have a space $H_\ast(\Phi(\{\lambda_1,\lambda_2\}))$. Any other object of
$\mscr{O}_\Lambda^d$ can be written as a union of objects of the form
$\{\lambda_1,\lambda_2\}$.  Since $\Phi$ is h-split, for any object $(O,s,t)$
of $\mscr{O}_\Lambda^d$, where $O$ is a non-negative integer,  and $s,t:  O \to
\Lambda$ are maps,
$$
H_\ast(\Phi(O,s,t)) = \otimes_{i =0}^{O-1} H_\ast(\Phi (
\{s(i),t(i)\}))
$$

Let $i : \mscr{O}_\Lambda^d \to \mscr{OC}_\Lambda^d$ and $j :
\mscr{C}^d \to \mscr{OC}_\Lambda^d$ denote the natural functors.
If $\Phi$ is an open-closed TCFT, then $j^\ast \Phi$ is a closed
TCFT and $i^\ast \Phi$ is an open TCFT.

Recall the objects of $\mscr{OC}_\Lambda^d$ are of the form $(C,O,s,t)$ where
$C,O$ are integers and $s,t: O \to \Lambda$ are maps. If $(\Lambda, \Phi)$ is
an open-closed TCFT, then $H_\ast(j^\ast \Phi)$ is the homology of the
associated closed TCFT, or equivalently the homology of $\Phi$ applied to the
object where $C = 1$  and $O  =0$.  Then,
$$
H_\ast(\Phi(C,O,s,t)) = \otimes_{o = 0}^{O - 1} H_\ast(\Phi(\{s(o),t(o)\})
\otimes H_\ast(j^\ast \Phi)^{\otimes C}
$$

Now we can state the main results of this paper.
\begin{thmA}
\begin{enumerate}
\item
The category of open TCFTs is quasi-equivalent to the category of
(unital) Calabi-Yau extended $A_\infty$ categories.
\item
Given any open TCFT, $(\Lambda,\Phi)$, we can push forward the
functor $\Phi : \mscr{O}^d_\Lambda \to \Comp$ to $\mbb{L} i_\ast
\Phi : \mscr{OC}^d_\Lambda \to \Comp$. This functor is again
h-split, so that $(\Lambda, \mbb{L} i_\ast \Phi)$ is an
open-closed TCFT.  This is the homotopy universal TCFT.
\item
We have a natural isomorphism
$$
H_\ast( j^\ast \mbb{L} i_\ast \Phi) \iso HH_\ast(A)
$$
where $A$ is the $A_\infty$ category corresponding to $(\Lambda,\Phi)$, and
$HH_\ast(A)$ is the Hochschild homology group.
\end{enumerate}
\end{thmA}
The notion of unital Calabi-Yau extended $A_\infty$ category will
be explained later.

\section{Combinatorial models for categories controlling open-closed topological conformal field theory}
\label{section models for moduli space}

In this section, an explicit dgsm category
$\mscr{D}^d_{\Lambda,open}$ is constructed which is quasi-isomorphic
to $\mscr{O}^d_\Lambda$.  This uses the cellular models for moduli
spaces which I introduced in \cite{cos_2004}, and which are
discussed in detail in \cite{cos_2006}.

The categories of modules for $\mscr{D}^d_{\Lambda,open}$ and
$\mscr{O}^d_\Lambda$ are quasi-equivalent. We have an $\op{SymOb}
\mscr{OC}_\Lambda -\mscr{O}^d_\Lambda$ bimodule,
$\mscr{OC}_\Lambda^d$. An explicit model $\mscr{D}^d_{\Lambda}$ for
the corresponding $\op{SymOb} \mscr{OC}_\Lambda
-\mscr{D}^d_{\Lambda,open}$ bimodule is constructed.

These results are enough to prove Theorem A. We will show later that a h-split
$\mscr{D}^{d}_{\Lambda,open}$ module is a unital extended Calabi-Yau $A_\infty$
category with set of objects $\Lambda$.  For each such, say $\Phi$, we will
calculate $H_\ast( \mscr{D}^d_{\Lambda} \otimes_{\mscr{D}^d_{\Lambda,open}}
\Phi)$, and find it is a tensor product of Hochschild homology groups of $\Phi$
and homology of morphism complexes of $\Phi$.  This will show that for any open
TCFT, the corresponding functor $\mscr{OC}^d_\Lambda \to \Comp$ is h-split and
has for homology of the closed states the Hochschild homology of the
corresponding $A_\infty$ category.

We do this by constructing cellular models for certain of our moduli
spaces of Riemann surfaces with open closed boundary.  Let $\alpha,
\beta \in \op{Ob} \mc M_{\Lambda}$ be such that $\alpha$ has no
closed part (so $\alpha = (0, O,s,t)$). We will construct
combinatorial models for the spaces $\mc M_{\Lambda}(\alpha,\beta)$.
The cell complex $\ocell(\alpha,\beta)$ we will construct will live
in a moduli space of Riemann surface with nodes along the boundary;
the surfaces in $\ocell(\alpha,\beta)$ will be those which are built
up from discs and annuli.

\subsection{A cellular model for moduli space}

The first step is to describe the moduli space of Riemann surfaces
with possibly nodal boundary. Let $\alpha,\beta \in \op{Ob} \mc
M_{\Lambda}$ as before, and assume $\alpha$ has no closed part.
Write $C(\beta),O(\beta),O(\alpha)$ for the closed and open
boundaries in $\alpha$ and $\beta$. Note that $C(\alpha) = 0$.

To keep the notation simple, I will omit the references to the
category of D-branes $\Lambda$, so that $\mc M$ will be synonymous
with $\mc M_\Lambda$.
\begin{definition}
Let $\pc (\alpha,\beta)$  be the moduli space of  Riemann surfaces $\Sigma$
with boundary,  with (outgoing) closed boundary components labelled by
$0,\ldots, C(\beta)-1$. These boundary components each have exactly one marked
point on them (this replaces the parameterisation on the boundary components of
the surfaces in $\oc$).  There are  further marked points labelled by
$O(\alpha)$ and $O(\beta)$ distributed along the remaining boundary components
of $\Sigma$; these correspond to the open boundaries. The free boundaries are
the intervals which lie between open boundaries, and the boundary components
with no marked points on them; these are labelled by D-branes in $\Lambda$ in a
way compatible with the maps $s,t : O(\alpha) \to \Lambda$ and $s,t : O(\beta)
\to \Lambda$.  Each connected component of the surface must have at least one
free boundary.

The surface $\Sigma$ may have nodes along the boundary, as in
\cite{liu2002}, \cite{cos_2004,cos_2006}. However, unlike in Liu's
work, there can be no nodes on the interior of $\Sigma$, nor are
there marked points on the interior. Marked points are not allowed
to collide with nodes.  Each closed boundary component of $\Sigma$
must be smooth (that is contain no nodes). Another difference from
Liu's work is that each boundary must be of positive length;
boundaries cannot shrink to punctures.   The surface $\Sigma$ must
be stable, that is have a finite automorphism group. This
corresponds to the requirement that no irreducible component of
$\Sigma$ can be a disc with $\le 2$ open marked points.

There are four exceptional kinds of surface; we allow surfaces
with connected components of this form.  The disc with zero, one
or two open marked points and the annulus with no open or closed
marked points are unstable; we declare the moduli space of any of
these types of surfaces to be a point.

\end{definition}

It is important to put in these  exceptional cases.  Part of $\pc$
will be made into a category, and the disc with one incoming and
one outgoing open point will be the identity.  The disc with one
open point will  give the unit in an $A_\infty$ category.

\begin{figure}
\includegraphics{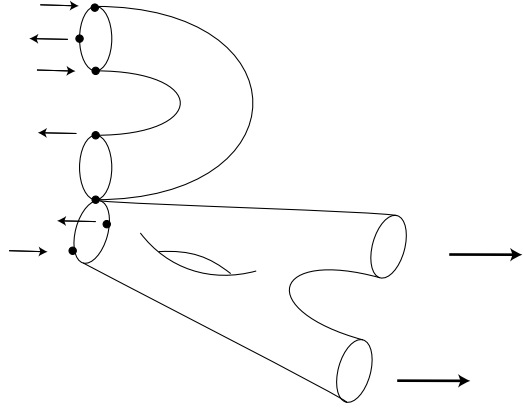}

\caption{A surface in $\pc (l)$.  The dots represent incoming or
outgoing open boundaries.  The boundaries with no dots are closed
and outgoing. }
\end{figure}

The moduli spaces $\pc$ are  orbifolds with corners. This follows from the work
of Liu \cite{liu2002}. One can see this by comparing the moduli spaces $\pc$ to
the real points of the Deligne-Mumford moduli spaces of curves. The interiors,
$\mc N$, are therefore smooth orbifolds. The spaces $\mc N$ parameterise
non-singular surfaces in $\pc$. The inclusion
$$\mc N \into \pc $$
is a $\Q$ homotopy equivalence.

The next step is to write down a subspace of the boundary of $\pc$
which is $\Q$ homotopy equivalent to $\pc$. Recall that the space
of (isomorphism classes of) annuli can be identified with $\mbb
R_{> 0}$.  Every annulus is isomorphic to an annulus of the form
$\{ z \mid 1 < \abs z < 1+R\}$ for some unique $R \in \R_{> 0}$,
which we call the modulus of the annulus.

\begin{definition}
Define $\ocell(\alpha,\beta) \subset \pc (\alpha,\beta)$ to be the subspace
consisting of surfaces $\Sigma \in \pc (\alpha,\beta)$, each of whose
irreducible components is either a disc, or an annulus of modulus $1$. We
require that one side of each annulus is an outgoing boundary component. Recall
that in $\pc (\alpha,\beta)$ the outgoing closed boundary components are
required to be smooth; this implies that the annuli are in one to one
correspondence with the outgoing closed boundary components $C(\beta)$.

$\ocell(\alpha,\beta)$ also contains the exceptional surfaces; we
allow surfaces with connected components which are discs with $\le
2$ marked points or annuli with no open or closed marked points.
\end{definition}

\begin{figure}

\includegraphics{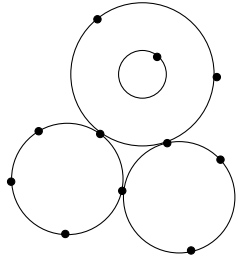}
\caption{A surface in $\ocell$, with $7$ open boundaries and one
closed boundary. The inside of the annulus is the outgoing closed
boundary component, the open boundaries may be incoming or
outgoing.}

\end{figure}

\begin{proposition}
The inclusion $\ocell(\alpha,\beta) \into \pc(\alpha,\beta)$ is a
weak homotopy equivalence of orbispaces (and therefore a $\Q$
homotopy equivalence of coarse moduli spaces).
\end{proposition}
\begin{proof}
This follows immediately from the results of \cite{cos_2004,
cos_2006}. For integers $g,h,r,s$ with $g,r,s \ge 0$, $h > 0$,
define the orbi-space $\pc_{g,h,r,s}$ to be the moduli space of
stable Riemann surfaces with possibly nodal boundary as above, with
$r$ boundary (open) marked points and $s$ internal marked points, of
genus $g$ with $h$ boundary components. As I discuss in detail in
\cite{cos_2006}, we have an orbi-cell complex $D_{g,h,r,s} \subset
\pc_{g,h,r,s}$ consisting of Riemann surfaces built up from discs,
each of which has at most one internal marked points. The inclusion
$D_{g,h,r,s} \into \pc_{g,h,r,s}$ is a weak homotopy equivalence.

We can replace the $s$ internal marked points by unparameterised
boundary components, in the moduli spaces $D_{g,h,r,s}$ and
$\pc_{g,h,r,s}$.  Evidently, all the corresponding moduli spaces are
homotopy equivalent, so the inclusion of these new spaces is also a
homotopy equivalence.. We can also add on to each of these $s$
boundary components a marked point, and the result continues to
hold, as we are simply passing to the total space of a torus bundle.

It follows immediately that the inclusion $G(\alpha,\beta) \into
\pc(\alpha,\beta)$ is a weak homotopy equivalence of orbispaces.

\end{proof}

Suppose $\alpha,\beta$ both satisfy $C(\alpha) = C(\beta) = 0$.
Then there are gluing maps
$$
\pc(\alpha,\beta) \times \pc(\beta,\gamma) \to \pc(\alpha,\gamma)
$$
These maps glue the outgoing open marked points of a surface in
$\pc(\alpha,\beta)$ to the corresponding incoming marked points of
a surface in $\pc(\beta,\gamma)$.  We need to describe how to glue
the exceptional surfaces; the discs with one or two marked points.
Gluing the disc with two open marked points, one incoming and one
outgoing, is the identity operation.  Gluing the disc with two
outgoing marked points onto two incoming marked points of a
surface $\Sigma$ corresponds to gluing the two marked points of
$\Sigma$ together; similarly for the disc with two incoming.
Gluing the disc with one marked point onto a marked point of a
surface $\Sigma$ causes us to forget that marked point.
\begin{lemma}
There is a category whose objects are the objects $\alpha$ of $\mc
M_\Lambda$ with $C(\alpha) = 0$ (i.e.\ no incoming closed
boundaries), whose morphisms are the spaces $\pc(\alpha,\beta)$ and
whose composition maps are the gluing described above.
\end{lemma}
Recall that $\mc M_\Lambda$ is the topological version of $\mscr{OC}_\Lambda$;
$\mscr{OC}_\Lambda$ is chains on $\mc M_\Lambda$. We defined  $\op{Ob} \mc
M_\Lambda$ to be the symmetric monoidal category with the same objects as $\mc
M_\Lambda$ but whose morphisms are the symmetry maps $a_1 \otimes \ldots a_n
\iso a_{\sigma(1)} \otimes \ldots \otimes a_{\sigma(n)}$, for $\sigma \in S_n$.

This lemma is clear. Call this category $\pc_{open}$. The spaces
$\ocell(\alpha,\beta) \subset \pc(\alpha,\beta)$ define a subcategory
$\ocell_{open} \subset \pc_{open}$. There is also the structure of symmetric
monoidal category on $\ocell_{open}$ and $\pc_{open}$ given by disjoint union.

$\pc$ defines a monoidal functor $\op{Ob} \mc M _\Lambda \times \pc_{open}^{op}
\to \op{Top}$, given by $(\beta,\alpha) \mapsto \pc(\alpha,\beta)$.  Similarly
$\ocell$ defines a functor $\op{Ob} \mc M_\Lambda \times \ocell^{op}_{open} \to
\op{Top}$.

 Let us take chain complexes
$C_\ast(\pc,\det^d)$ where $\det^d$ is the local system defined
before; $C_\ast(\pc_{open},\det^d)$ is a differential graded
symmetric monoidal category, and $C_\ast(\pc,\det^d)$ defines an
$\op{Ob} \mscr{OC}^d - C_\ast(\pc_{open},\det^d)$ bimodule.
\begin{proposition}
The dgsm category $C_\ast(\pc_{open},\det^d)$ is quasi-isomorphic
to the dgsm category $\mscr{O}^d$.

Under the induced quasi-equivalence of categories between $\op{Ob} \mscr{OC}^d
- C_\ast(\pc_{open},\det^d)$ bimodules and $\op{Ob} \mscr{OC}^d - \mscr{O}^d$
bimodules, $C_\ast(\pc,\det^d)$ corresponds to $\mscr{OC}^d$.
\end{proposition}
We are suppressing the set $\Lambda$ of D-branes from the notation
here.
\begin{proof}
I will sketch the proof of the statement about categories, in the
case $d = 0$; the remaining statements are proved in a similar
way.  We will do the topological version, and find a topological
category $\til {\oc}_{open}$, with the same objects as
$\oc_{open}$, and with functors $\pc_{open} \to \til{ \mc
M}_{open} \from \mc M_{open}$ which are $\Q$ homotopy equivalences
on the spaces of morphisms.

For $\alpha,\beta \in \op{Ob} \mc M$, let $\til{ \mc
M}_{open}(\alpha,\beta)$ be the moduli space of surfaces with
nodal boundary, as in $\pc(\alpha,\beta)$, but now the open
boundaries are parameterised embedded intervals, like in
$\oc(\alpha,\beta)$.  These intervals do not intersect the nodes
or each other.  Each outgoing open boundary has a number $t \in
[0,1/2]$ attached to it.

The gluing which defines the maps $\til{ \mc
M}_{open}(\alpha,\beta) \times \til{ \mc M}_{open}(\beta,\gamma)
\to \til{ \mc M}_{open}(\alpha,\gamma)$ is defined as follows. Let
$\Sigma_1 \in \til{ \mc M}_{open}(\alpha,\beta)$ and $\Sigma_2 \in
\til{ \mc M}_{open}(\beta,\gamma)$, and let $o \in O(\beta)$. This
corresponds to an open boundary on each of of the $\Sigma_i$. Let
$t \in [0,1/2]$ be number corresponding to $o$.  Glue the
subinterval $[t,1-t] \subset [0,1]$ of the boundary on $\Sigma_1$
to the corresponding subinterval $[t,1-t]$ of the corresponding
boundary on $\Sigma_2$.

This evidently makes $\til{ \mc M}_{open}$ into a category.  The map
$\mc M_{open} \into \til{ \mc M}_{open}$ assigns the number $0$ to
the open boundaries, and is a homotopy equivalence on spaces of
morphisms.  Similarly, the map $\pc_{open} \into \til{ \mc
M}_{open}$ assigns the number $1/2$ to open boundaries, and is a
homotopy equivalence on the space of morphisms.

This argument implies the corresponding result at chain level, and
extends without difficulty to the case of twisted coefficients and
to yield an equivalence of modules.

\end{proof}

We want to give an orbi-cell decomposition of the spaces
$\ocell$.  We will do this by writing down a stratification of
$\ocell$ whose strata are orbi-cells, that is the quotient of a
cell by a finite group.  There is an obvious stratification of
$\ocell$, given by the topological isomorphism type of the
corresponding marked nodal surface.  This is not quite a cell
decomposition, as the moduli space of marked points on the
boundary of an annulus, one of whose boundaries is closed, is not
contractible, but is homotopic to $S^1$. We need to refine this
stratification a little.

Let $\Sigma \in \ocell(\alpha,\beta)$.  Let us assume for
simplicity that no connected component of $\Sigma$ is an
exceptional (unstable) surface. We will give $\Sigma$ a cell
decomposition. Let $A \subset \Sigma$ be an irreducible component
which is an annulus with a closed boundary. In order to get a cell
decomposition on $\Sigma$, we have to make a cut on the annulus.
Let $A_{closed}, A_{open}$ be the boundary components of $A$;
where $A_{closed}$ has precisely one marked point, $p$ say,
corresponding to an outgoing closed boundary of $\Sigma$, and
$A_{open}$ may have several incoming and outgoing open marked
points and possibly some nodes. There is a unique holomorphic
isomorphism from $A$ to the cylinder $S^1 \times [0,1]$,  such
that  $p \in A_{closed}$ goes to $(1,0) \in S^1 \times [0,1]$. The
inverse image of $1 \times [0,1]$ in this gives a cut on the
annulus, starting at $p \in A_{closed}$ and ending at some point
$p'$ on $A_{open}$. Now give $\Sigma$ a cell decomposition, by
declaring that the $0$ skeleton consists of the nodes, marked
points, and the places where the cut on an annulus intersect the
boundary of the annulus; the one cells are $\partial \Sigma$,
together with the cuts on the annuli; and the $2$ skeleton is
$\Sigma$.   The two cells of $\Sigma$ are oriented, and $\Sigma$
is marked by D-branes, incoming/outgoing open marked points, and
closed outgoing marked points.

\begin{figure}

\subfigure[] {
\includegraphics{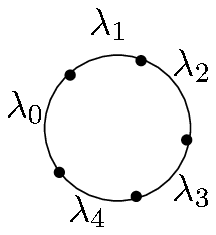}
\label{figure disc cell} }  \subfigure[]{
\includegraphics{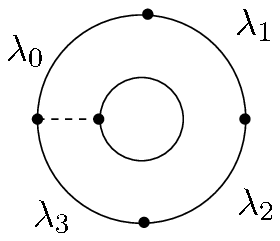} \label{figure annulus cell1} }
\subfigure[]{ \includegraphics{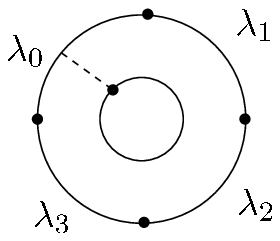} \label{figure annulus
cell2}} \centering \caption{}
\begin{minipage}[b]{120mm}
\vspace{3mm} These are the three basic types of cell in the moduli
spaces $\ocell$, from which all others are built by open gluing.

\vspace{4mm} Figure \ref{figure disc cell} represents the cell in
moduli space of points moving on a disc. The marked points are
open boundaries, and may be incoming or outgoing;  the $\lambda_i$
are D-brane labellings on free boundaries.

\vspace{4mm} Figures \ref{figure annulus cell1} and \ref{figure annulus cell2}
are the two kinds of cell in the space of marked points on the annulus. The
interior of the annulus is a closed outgoing boundary; the marked point on this
represents the start of the closed boundary. The remaining marked points are
open, incoming or outgoing.

\vspace{4mm} In figure \ref{figure annulus cell1}  the closed marked point is
parallel to an open one, where as in figure \ref{figure annulus cell2}, the
closed marked point is parallel to the interior of a free boundary.
\end{minipage}
\end{figure}

Give $\ocell(\alpha,\beta)$ a stratification by saying that two
surface $\Sigma_1, \Sigma_2$ are in the same stratum if and only
if the corresponding marked, oriented 2-cell complexes in
$\adm(\alpha,\beta)$ are isomorphic.
\begin{lemma}
This stratification of $\ocell(\alpha,\beta)$ is an orbi-cell
decomposition, and further the composition maps
$\ocell(\alpha,\beta)\times\ocell(\beta,\gamma) \to
\ocell(\alpha,\gamma)$ are cellular. \label{lemma cell}
\end{lemma}
To show that this stratification is an orbi-cell decomposition,
the main point to observe is that the stratification of the space
of marked points on the annulus is indeed a cell decomposition.

We are using a strong notion of cellular map : a map $f : X \to Y$ between
(orbi)-cell complexes is cellular if $f^{-1} Y_i = X_i$, where $X_i$ is the
union of cells of dimension $\le i$.

Define
$$\choc(\alpha,\beta) = C_\ast^{cell} (\ocell(\alpha,\beta)) \otimes \K$$
to be the associated complex of $\K$ cellular chains. Similarly,
for an integer $d \ge 0$, define
$$
\choc^d(\alpha,\beta)  = C_\ast^{cell} (\ocell(\alpha,\beta),
\op{det}^d ) \otimes \K
$$
Here we take cellular chains with local coefficients.

Let us describe informally the chain complexes
$\choc(\alpha,\beta)$.   Each Riemann surface in
$\ocell(\alpha,\beta)$ determines a cell in the moduli space, and
so an element of the cellular chain group $\ocell(\alpha,\beta)$.
Thus, we can think of a chain in  $\choc(\alpha,\beta)$ as being
represented by a surface, and similarly for
$\mscr{D}^d(\alpha,\beta)$. The boundary maps in
$\choc^d(\alpha,\beta)$ correspond to degenerating surfaces to
allow more nodes, and also allowing a closed marked point, on the
boundary of an annulus, to become parallel to an open marked point
or node on the other boundary of an annulus.

There are composition maps $\choc^d(\alpha,\beta) \otimes \choc^d(\beta,\gamma)
\to \choc^d(\alpha,\gamma)$, which make $\choc^d_{open}$ (the part where
$\alpha,\beta$ have only open boundaries) into a differential graded symmetric
monoidal category, and $\choc^d$ into a $\op{Ob} \mscr{OC}^d - \choc^d_{open}$
bimodule.
\begin{lemma}
The differential graded symmetric monoidal categories $\choc^d_{open}$ and
$\mscr{O}^d$ are quasi-isomorphic. Under the induced quasi-equivalence of
categories between  \linebreak$\op{Ob} \mscr{OC}^d - \choc^d_{open}$ bimodules
and $\op{Ob} \mscr{OC}^d - \mscr{O}^d$ bimodules, $\choc^d$ corresponds to
$\mscr{OC}^d$.
\end{lemma}
The point is that the chain complex functor $C_\ast$ constructed
in the appendix has the property that for each (orbi)-cell complex
$X$, there is a quasi-isomorphism $C_\ast^{cell}(X) \to
C_\ast(X)$, compatible with products and natural for cellular
maps.  The same holds when we take chains with local coefficients.
This shows that the functor $\choc^d_{open} \to C_\ast
(\ocell_{open},\det^d)$ is a quasi-isomorphism, and we have
already seen that $C_\ast(\ocell,\det^d) \simeq
C_\ast(\pc_{open},\det^d) \simeq \mscr{O}^d$.  Similar remarks
prove the statement about $\choc^d$ as an  $\op{Ob} \mscr{OC}^d -
\choc^d_{open}$ bimodule.

\subsection{Generators and relations for  $\choc^d_{open}$}

If $\lambda_0,\ldots,\lambda_n$ is an ordered set of D-branes, let
$\{\lambda_0,\ldots,\lambda_n \}$ be the object in $\op{Ob} \choc^d_{open}$
with $O = n$ and $s(i) = \lambda_i$, $t(i) = \lambda_{i+1}$ for $0 \le i \le
O-1$.

Use the notation
$$
\{\lambda_0,\ldots,\lambda_{n-1}\}^c =
\{\lambda_0,\ldots,\lambda_{n-1},\lambda_0\}
$$
The superscript $c$ stands for cyclic.

Define an element $D(\lambda_0,\ldots,\lambda_{n-1})$ of
$\choc^d_{open}(\{\lambda_0,\ldots,\lambda_{n-1}\}^c,0)$, given by
the cellular chain which is the disc with $n$ marked points on it,
all incoming, with the cyclic order $0,1,\ldots,n-1$, labelled in
the obvious way by D-branes; as in figure \ref{figure disc cell}.
(Pick, arbitrarily, some orientation on this cell, and a section
of $\det^d$, in order to get a  cellular chain).

Note that $D(\lambda_0, \ldots,\lambda_{n-1})$ is cyclically
symmetric up to sign; so that
$$D(\lambda_0, \ldots,\lambda_{n-1})
= \pm D(\lambda_1, \ldots,\lambda_{n-1}, \lambda_0)$$ under the
permutation isomorphism between
$\{\lambda_0,\ldots,\lambda_{n-1}\}^c$ and
$\{\lambda_1,\ldots,\lambda_{n-1},\lambda_0\}^c$.

When $n \ge 3$, $D(\lambda_0,\ldots,\lambda_{n-1})$ is an element of degree $n
- 3 + d$. When $n = 1,2$ it is an element of degree $d$.

Let $C \subset \choc^d_{open}$ be the subcategory with the same objects,  but
whose morphism surfaces are not allowed to have connected components which are
the disc with $\le 1$ open marked points; or the disc with two open marked
points, both incoming; or the annulus with neither open or closed marked
points.  We consider the morphisms in $C$ not to be complexes, but to be graded
vector spaces : we forget the differential.
\begin{proposition}
$C$ is freely generated, as a symmetric monoidal category over the
symmetric monoidal category $\op{Ob} \choc^d_{open}$, by the discs
$D(\lambda_0,\ldots,\lambda_{n-1})$, where $n \ge 3$, and the
discs with two outgoing marked points, subject to the relation
that $D(\lambda_0,\ldots,\lambda_{n-1})$ is cyclically symmetric
(up to an appropriate sign). \label{prop free generators}
\end{proposition}
The sign in the cyclic symmetry is determined by the choice of
orientation on the cell in $\ocell$ corresponding to
$D(\lambda_0,\ldots,\lambda_{n-1})$.

Note that it makes sense to talk about generators and relations for a symmetric
monoidal category; this is because we have fixed the base category $\op{Ob}
\choc^d_{open}$, and the new symmetric monoidal category we are constructing
has the same set of objects.  The morphism spaces of a symmetric monoidal
category given by generators and relations will be built up using composition,
tensor product, and adding morphisms of the base category $\op{Ob}
\choc^d_{open}$,  from the generators.

\begin{proof}[Proof of proposition \ref{prop free generators}]
Let $C'$ be the category with these generators and relations.  There is a
functor $C' \to C$; firstly we will show this is full.  We can take disjoint
union of surfaces in $C'$, and we can use the disc with two outgoing marked
points to change an incoming boundary of a surface to an outgoing boundary.
Clearly, any surface in $C(\alpha,\beta)$ can be built up using disjoint union
and gluing from discs. This shows $C' \to C$ is full.

Next, we need to show that this functor is faithful. It suffices to write down
a functor $C \to C'$ which is an inverse. On objects, this is the identity.
Suppose we have a surface $\Sigma$ in $C(\alpha,\beta)$. We can write
$$\Sigma =  \Sigma' \circ \phi$$
in a unique way, where $\Sigma'$ is a disjoint union of identity maps and discs
with all incoming boundaries, and $\phi$ is a disjoint union of discs with two
outgoing boundaries and identity maps.  $\Sigma'$ is the normalisation of
$\Sigma$ with all of its marked points made incoming. $\phi$ has the effect of
gluing the marked points of $\Sigma'$ which correspond to nodes of $\Sigma$
together, and of changing the incoming points of $\Sigma'$ which correspond to
outgoing points of $\Sigma$ into outgoing.

This decomposition of $\Sigma$ allows us to write down the inverse map
$C(\alpha,\beta) \to C'(\alpha,\beta)$, and it is easy to check this defines a
functor.
\end{proof}

Let $\choc^+_{open} \subset \choc^d_{open}$ be the subcategory
with the same objects but whose morphisms are given by disjoint
unions of discs, with each connected component having precisely
one outgoing boundary.  Note that this is indeed a subcategory,
and is also independent of $d$; the local system $\det^d$ can be
canonically trivialised in degree $0$ on the moduli space of discs
with one outgoing boundary.

\begin{figure}

\includegraphics{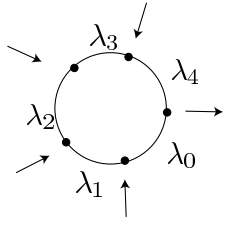}

 \caption{The chain $D^+(
 \lambda_0,\lambda_1,\lambda_2,\lambda_3,\lambda_4)$.} \label{figure disc tree level}
\end{figure}

For each ordered set $\lambda_0,\ldots, \lambda_{n-1}$ of
D-branes, where $n \ge 1$, let $D^+(\lambda_0,\ldots,
\lambda_{n-1})$ be the disc with $n$ marked points, and D-brane
labels by the $\lambda_i$, but such that all of the marked points
are incoming except that between $\lambda_{n-1}$ and $\lambda_0$,
as in figure \ref{figure disc tree level}.

$D^+(\lambda_0,\ldots,\lambda_{n-1})$ is in $\op{Hom} (
\{\lambda_0,\ldots,\lambda_{n-1} \} , \{ \lambda_0,\lambda_1 \})$.

\begin{lemma}
$\choc^+_{open}$ is freely generated, as a symmetric monoidal
category over $\op{Ob} \choc^d_{open}$, by the discs
$D^+(\lambda_0,\ldots,\lambda_{n-1})$, modulo the relation that
$$
 D^+ ( \lambda_0, \ldots,\lambda_i, \lambda_i,
\ldots, \lambda_{n-1} )\circ  D^+(\lambda_i) = 0
$$
whenever $n \ge 4$, and when $n = 3$,
\begin{align*}
 D^+ ( \lambda_0,\lambda_0,\lambda_1) & \circ D^+ ( \lambda_0 ) \\
 D^+(\lambda_0,\lambda_1,\lambda_1) & \circ D^+(\lambda_1)
\end{align*}
are both the identity map on the object $\{\lambda_0,\lambda_1\}$.
\label{lemma generators relations operad}
\end{lemma}
This is basically a corollary of the previous result. Note that the relations
stated do indeed hold; composing with $D^+(\lambda_i)$ has the effect of
forgetting the open marked point which lies between the two copies of
$\lambda_i$.  By ``composing'' we mean of course placing the identity on all
other factors.

\begin{theorem}
$\choc^d_{open}$ is freely generated, as a symmetric monoidal
category over $\op{Ob} \choc^d_{open}$, by $\choc^+_{open}$, and
the discs with two incoming or two outgoing boundaries, modulo the
following relations.

The first relation is illustrated in the figure \ref{figure ip
relation} ; it says that an appropriate gluing of the disc with
two outgoing boundaries and with two incoming boundaries yields
the identity (a disc with one incoming and one outgoing boundary).

Observe that we can change an outgoing boundary to an incoming
boundary; let $D(\lambda_0,\ldots,\lambda_{n-1})$ be obtained from
$D^+(\lambda_0,\ldots,\lambda_{n-1})$ by making the outgoing
boundary incoming.  The second relation is that
$D(\lambda_0,\ldots,\lambda_{n-1})$ is cyclically symmetric (up to
an appropriate sign).

\label{theorem generators relations category}
\end{theorem}
\begin{figure}

\includegraphics{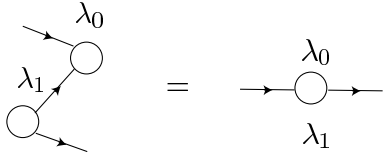}
\caption{ Gluing of a disc with two outgoing to a disc with two
incoming yields the identity, a disc with an incoming and and
outgoing. }

\label{figure ip relation}
\end{figure}

This follows almost immediately from the previous result. This generators and
relations description of course refers to the category without the
differential.  Note that the disc with no marked points and the annulus with no
open or closed marked points are included in $\choc^d_{open}$; for example, the
annulus with no marked points is given by gluing the disc with two outgoing
marked points to that with two incoming marked points.

Let $\lambda_0,\ldots,\lambda_{n-1}$ be an ordered set of
D-branes.  There is an element
$$A(\lambda_0,\ldots,\lambda_{n-1}) \in  \choc^d(\{\lambda_0,\ldots,\lambda_{n-1}\}^c, (1,0) )$$ given by
given by the annulus with $n$ marked points, and the intervals
between the marked points labelled by the D-branes $\lambda_i$, as
in figure \ref{figure annulus cell1}. The parameterisation on the
closed boundary - on the interior of the annulus - starts at the
open marked point $0$ between $\lambda_{n-1}$ and $\lambda_0$.

\begin{figure}

\includegraphics{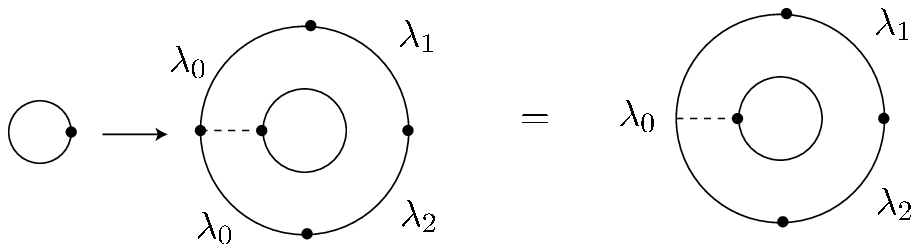}
 \caption{} \label{figure annulus forget marked point}
\end{figure}

The object $(1,0)$ of $\mscr{OC}^d$ has one closed boundary and no
open boundaries.   Note that $A(\lambda_0,\ldots,\lambda_{n-1})$
is an $n-1$ chain in $\choc^d$.
\begin{theorem}
The $\op{Ob} \mscr{OC} - \choc^d_{open}$ bimodule $\choc^d$ is freely
generated,  by the $A(\lambda_0,\ldots,\lambda_{n-1})$, and the identity maps
$1 \in \choc^d_{open}(\alpha,\alpha) \subset \choc^d(\alpha,\alpha)$, modulo
the following relations.  Firstly, if we glue the disc with one boundary to any
of the open marked points of $A(\lambda_0,\ldots,\lambda_{n-1})$, except that
lying between $\lambda_{n-1}$ and $\lambda_0$, we get $0$. Secondly, the
disjoint union of the identity element on $\alpha$ with that on $\beta$ is the
identity  on $\alpha \amalg \beta$.

\label{theorem open closed generators relations}
\end{theorem}
This is proved in essentially the same way that the previous results are.  The
main points are as follows. Since $\choc^d$ is an $\op{Ob} \mscr{OC} -
\choc^d_{open}$ bimodule, we can take disjoint unions, so we get disjoint
unions of annuli and identity elements.  We also get discs using the action of
$\choc^d_{open}$ on the identity elements. For example, we have the identity
element for the zero object $\alpha = 0$, which gives us discs with all
incoming boundaries.  The action of $\choc^d_{open}$ allows us to glue discs to
annuli.  This also lets us glue annuli together,  and change incoming marked
points on annuli to outgoing, using the disc with two outgoing boundaries. One
point to observe is that if we glue a disc with one marked point to the annulus
$A(\lambda_0,\ldots,\lambda_{n-2},\lambda_{0})$ at the marked point between
$\lambda_{n-1}$ and $\lambda_0$ we get an annulus where the starting point for
the parameterisation of the closed boundary lies in the free boundary
$\lambda_0$, as in figure \ref{figure annulus forget marked point}.

 This ensures that although the moduli space of annuli
contains two types of cells, depending on whether the start of the
parameterisation on the closed boundary is at an open or a free
boundary, we need only take one type as a generator.

\begin{definition}
Let $\choc^+$ be the $ \op{Ob} \mscr{OC} - \choc^+_{open} $
bimodule with the same generators and relations as $\choc^d$.
\label{definition dplus}
\end{definition}
Note that this makes sense, as the relations involve only the disc
with one outgoing marked point, which comes from $\choc^+_{open}$.
It is clear that
$$
\choc^d = \choc^+ \otimes_{\choc^+_{open}} \choc^d_{open}
$$
as a $\op{Ob} \mscr{OC} - \choc^d_{open}$ bimodule.  Further, for
any left $\choc^d_{open}$ module $M$,
$$
\choc^d \otimes_{\choc^d_{open}} M = \choc^+
\otimes_{\choc^+_{open}} \choc^d_{open} \otimes_{\choc^d_{open}} M
= \choc^+ \otimes_{\choc^+_{open}} M
$$

\subsection{The differential in $\choc^d$}

We also want to describe the differential in the complexes
$\choc^d$. This is characterised by the fact that it respects the
composition maps $\choc^d(\alpha,\beta) \otimes
\choc^d(\beta,\gamma) \to \choc^d(\alpha,\gamma)$, and the way it
behaves on discs and annuli, which are the generators.   I will
only write down the formula up to sign; the precise signs will
depend on the orientation chosen for the cells in $\ocell$ of
marked points on discs and annuli. The precise signs don't matter.

The differential on discs is shown in figure \ref{figure
differential discs}.
\begin{figure}

\includegraphics{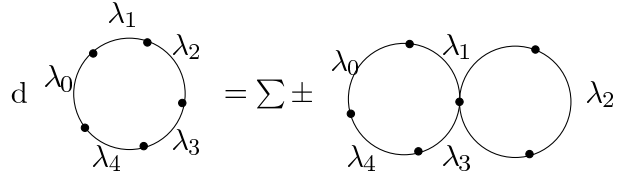}
\caption{The differential of a chain given by marked points on a
disc.  The marked points may be incoming or outgoing; the
$\lambda_i$ are D-branes.} \label{figure differential discs}
\end{figure}
This can be written as
$$
\d D(\lambda_0,\ldots,\lambda_{n-1}) = \sum_{\substack {0 \le i
\le  j \le n-1 \\ j - i \ge 2} } \pm D(\lambda_i,\ldots ,
\lambda_j ) \ast D(\lambda_j, \ldots, \lambda_i)
$$
where the $\ast$ indicates that we glue the open marked points
between $\lambda_i$ and $\lambda_j$ on each disc together.

\begin{figure}

\includegraphics{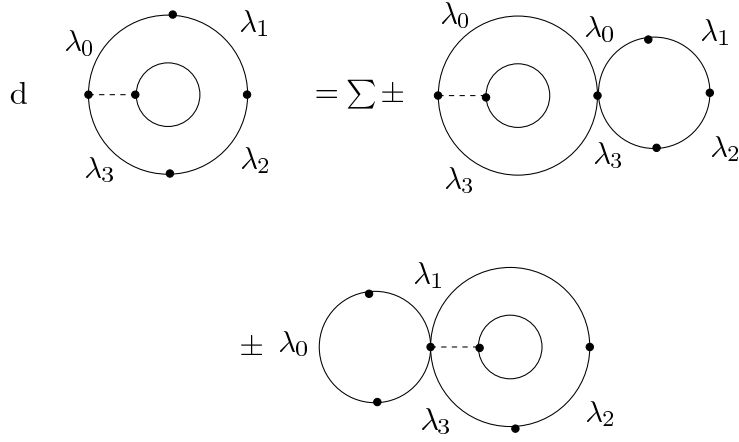}
 \caption{The differential of a
chain given by marked points on an annulus.  The interior circle
of the annulus is a closed outgoing boundary, the marked points on
the exterior may be incoming or outgoing open, and the $\lambda_i$
are D-branes.} \label{figure differential annuli}
\end{figure}

On annuli, it is given in figure \ref{figure differential annuli}.
This can be written as
\begin{multline}
\d A ( \lambda_0,\ldots, \lambda_{n-1})  = \sum_{ \substack { 0
\le i < j \le n-1  \\ \abs{i - j} \ge 2}  } \pm A( \lambda_0,
\ldots, \lambda_{i-1},\lambda_i, \lambda_j, \lambda_{j+1},
\ldots, \lambda_{n-1} ) \ast D(\lambda_i, \ldots, \lambda_j ) \\
 + \sum_{\substack { 0 \le j \le i \le n-1 \\ (j,i) \neq (0,n-1) } }  \pm A(\lambda_{j},\ldots,
 \lambda_i ) \ast D ( \lambda_i, \ldots, 0, 1, \ldots, \lambda_j )
\end{multline}
where, as before, the symbol $\ast$ means we should glue at the
open marked points between the D-branes $\lambda_i$ and
$\lambda_j$.

\begin{lemma}
\begin{enumerate}
\item
The $\op{Ob} \mscr{OC}^d-\choc^d_{open}$ bimodule $\choc^d$ is
$\choc^d_{open}$-flat.
\item
If $M$ is a h-split $\choc^d_{open}$ module, then
$$
\choc^d \otimes_{\choc^d_{open}}  M
$$
is a h-split $\op{Ob} \mscr{OC}^d$-module.
\end{enumerate}
The same is true if we consider $\choc^+_{open}$ and $\choc^+$
instead of $\choc^d_{open}$ and $\choc^d$. \label{lemma flat and
split}
\end{lemma}
\begin{proof}
Recall $\choc^d$ is generated as a $\choc^d_{open}-\op{Ob}
\mscr{OC}^d$ bimodule by the annuli
$A(\lambda_0,\ldots,A_{\lambda_{n-1}})$, and the identity elements
in $\choc^d(\alpha,\alpha)$ where $\alpha \in \op{Ob}
\choc^d_{open}$.

Filter $\choc^d$ as a bimodule, by giving a filtration on the generators,
defined by saying each identity element in $\choc^d(\alpha,\alpha)$ is in $F^0$
and each annulus $A(\lambda_0,\ldots,\lambda_{n-1})$ is in $F^n$.  The formula
for the differential of the annuli guarantees that this is a filtration as
complexes; indeed, $\d A(\lambda_0,\ldots,\lambda_{n-1})$ is in $F^{n-1}$.

Let $M$ be a left $\choc^d_{open}$ module.  Suppose $M \to M'$ is
a quasi-isomorphism; we want to show that the map
\begin{equation}
\choc^d (-,\beta) \otimes_{\choc^d_{open}}  M \to \choc^d
(-,\beta) \otimes_{\choc^d_{open}}  M' \label{equation qis}
\end{equation}
 is a quasi-isomorphism. Give
both sides the filtration induced from that on $ \choc^d
(-,\beta)$; it suffices to show that the map on the associated
graded complexes is a quasi-isomorphism.

This follows immediately from the generators and relations
description of $\choc^d$.  Let $\alpha \in \op{Ob}
\choc^d_{open}$; for an integer $C$, $C \amalg \alpha \in \op{Ob}
\mscr{OC}^d$; we add on $C$ closed states.  We want to show that
the map \ref{equation qis} is a quasi-isomorphism, with $\beta = C
\amalg \alpha$.  For simplicity I will show this when $C = 1$.

Then,  $\op{Gr}^n M \otimes_{\choc^d_{open}} \choc^d ( -, \alpha
\amalg 1)$ is spanned by the spaces
$$
M ( \alpha \amalg \{\lambda_0,\ldots,\lambda_{n-1}\}^c)
$$
This corresponds to putting the generators of $\choc^d$ which are
the identity in $\choc^d_{open}$ on the $\alpha$ factor and the
annulus $A(\lambda_0,\ldots,\lambda_{n-1})$ on the
$\{\lambda_0,\ldots,\lambda_{n-1}\}^c$ factor.

The only relation is that the composed map
\begin{multline}
M ( \alpha \amalg \{\lambda_0,\ldots,\lambda_{i-1}, \lambda_{i+1}
, \lambda_{n-1}\}^c ) \to M ( \alpha \amalg
\{\lambda_0,\ldots,\lambda_{i-1}, \lambda_i, \lambda_i,
\lambda_{i+1} , \lambda_{n-1}\}^c ) \\
\to \op{Gr}^n  \choc^d ( -,  \alpha \amalg 1)
\otimes_{\choc^d_{open}} M \label{equation annuli relation}
\end{multline}
is zero.  The first map comes from the element of
$$
\choc^d_{open} (\alpha \amalg \{\lambda_0,\ldots,\lambda_{i-1},
\lambda_{i+1} , \lambda_{n-1}\}^c ,\alpha \amalg
\{\lambda_0,\ldots,\lambda_{i-1}, \lambda_i, \lambda_i,
\lambda_{i+1} , \lambda_{n-1}\}^c  )
$$
which is the tensor product of the identity on $\alpha$ and
${\lambda_0,\ldots,\lambda_{i-1}, \lambda_{i+1} , \lambda_{n-1}}^c$ and the map
$0 \to \{\lambda_i,\lambda_i\}$ given by the disc with one outgoing marked
point.

The first map in the diagram  \ref{equation annuli relation} is always
injective; we can find a splitting coming from the disc with one incoming
marked point.  Thus the operation of taking the quotient is exact.

There is a similar description of $\choc^d(-,\beta)
\otimes_{\choc^d} M$, for all $\beta$, and the same argument shows
that the functor $\choc^d(-,\beta)\otimes_{\choc^d_{open}} -$ is
exact.  This proves the first part of the lemma.

The second part of the lemma is proved in a similar way. Let
$$
N = \choc^d \otimes_{\choc^d_{open}} M
$$
Then the filtration on $\choc^d$ induces one on $N$.  To show the
maps $N(\beta)\otimes N(\beta') \to N(\beta \amalg \beta')$ are
quasi-isomorphism, it suffices to do so on the associated graded.
This follows immediately from the description of $N$ given above.

Exactly the same proof shows the corresponding results for
$\choc^+_{open}$ and $\choc^+$.
\end{proof}

\section{Proof of the main results}

\label{section proof}
\subsection{$A_\infty$ categories}

 Let us recall some details of the definition of an $A_\infty$ category $D$.  There is a set
$\op{Ob} D$ of objects, and for each pair $A,B$ of objects, a
finite dimensional complex of $\K$ vector spaces $\op{Hom}(A,B)$.
The homological grading convention is used, so that the
differential is of degree $-1$.  For each sequence $A_0, \ldots ,
A_{n}$ of objects, where $n \ge 2$, there are maps
$$
m_n : \op{Hom} (A_0, A_1) \otimes \ldots \otimes \op{Hom}(A_{n-1},
A_{n}) \to \op{Hom}(A_0, A_{n})
$$
of degree $n-2$.  (Note this is different from the standard
convention of $2-n$). The differential on the complex
$\op{Hom}(A,B)$ is $m_1$. These maps must satisfy identities of
the form
\begin{multline*}
\sum_{0 \le i \le j \le n-1} \pm  m_{n - j + i} \left( \alpha_0
\otimes \ldots \otimes \alpha_{i-1} \otimes m_{j - i + 1} (
\alpha_{i} \otimes \ldots \otimes \alpha_{ j} ) \otimes
\alpha_{j+1 } \otimes \ldots \otimes \alpha_{n-1} \right) = 0
\end{multline*}

All our $A_\infty$ categories will be unital.  A unital $A_\infty$
category is an $A_\infty$ category $D$, together with for each $A
\in \op{Ob} D$  a closed element $1_A \in \op{Hom}_0(A,A)$, with
the following properties. Firstly,
\begin{align*}
m_2 ( \alpha \otimes 1_A ) &= \alpha \\
m_2 ( 1_A \otimes \beta )&= \beta
\end{align*}
for any $\alpha : B \to A$ and $\beta : A \to B$. Secondly, if
$\alpha_i : A_i \to A_{i + 1 }$ are maps,  for $0 \le i < n$, and
if $j = j+1$, then
$$
m_n ( \alpha_0 \otimes \alpha_1 \otimes \ldots 1_{A_j} \otimes
\ldots \alpha_{n-1} ) = 0
$$

\subsection{Calabi-Yau $A_\infty$ categories} A \emph{Calabi-Yau}
$A_\infty$ category of dimension $d$ is an $A_\infty$ category
$D$, with for each pair $A,B \in \op{Ob} D$, a closed
non-degenerate pairing
$$
\ip{ \quad}_{ A,B} : \op{Hom}(A,B) \otimes \op{Hom}(B,A) \to \K[d]
$$
which is symmetric (in the sense that $\ip{}_{A,B} = \ip{}_{B,A}$
under the natural symmetry isomorphism $\op{Hom}(A,B) \otimes
\op{Hom}(B,A) \iso \op{Hom}(B,A)\otimes \op{Hom}(A,B)$), and such
that the cyclic symmetry identity identity
$$
\ip { m_{n-1} ( \alpha_0 \otimes \ldots \otimes \alpha_{n-2}),
\alpha_{n-1} }  = (-1)^{(n+1) + \abs{\alpha_0} \sum_{i = 1}^{n-1}
\abs {\alpha_i}} \ip { m_{n-1} ( \alpha_1 \otimes \ldots \otimes
\alpha_{n-2}), \alpha_{0} }
$$
holds.

Let $X$ be a smooth projective Calabi-Yau variety of dimension
$d$. Let $\mc D^b(X)$ be the bounded derived category of coherent
sheaves on $X$. Then $\mc D^b(X)$ is a unital Calabi-Yau
$A_\infty$ category, of dimension $d$.  We have to change the
grading, so that
$$
\op{Hom}_i (A,B) = \op{Ext}^{-i}(A,B)
$$
The composition maps $m_n$ are all zero for $n \neq 2$.  The
pairing
$$\op{Hom}_{i} (A,B)^{\vee} \iso \op{Hom}_{-d -i } (B,A)$$
is Serre duality. (We need to trivialise the Serre functor, by
picking a non-zero holomorphic top form.).

This should not be regarded as being the correct Calabi-Yau
$A_\infty$ category for the $B$-model, as I mentioned in the
introduction. We should use an $A_\infty$ version $\mc
D^b_\infty(X)$.

\subsection{Open topological conformal field theories and $A_\infty$ categories}
Let $\Lambda$ be a set of D-branes.

Recall a monoidal functor between monoidal categories is called
\emph{split} if the maps
$$
F(a) \otimes F(b) \to F(a \otimes b)
$$
are isomorphisms.
\begin{lemma}
A split functor $\Phi : \choc^+_{open,\Lambda} \to \Comp$ is the same as a
unital $A_\infty$ category with set of objects $\Lambda$.
\end{lemma}
\begin{proof}
 Let $\Phi : \choc^+_{open,\Lambda} \to
\Comp$ be a split symmetric monoidal functor.  Then for each integer $O$, with
D-brane labels $s(i),t(i)$, for $0 \le i \le O-1$, we have a
natural isomorphism
$$
\Phi ( O,s,t) \iso \otimes_{i = 0}^{O-1} \Phi ( \{s(i), t(i)\})
$$
For each pair $\lambda,\lambda'$ of D-branes, write $\op{Hom}
(\lambda,\lambda') = \Phi ( \{\lambda,\lambda'\})$.

Generators and relations for  $\choc^+_{open,\lambda}$ are given
in \ref{lemma generators relations operad}.  The discs
$D^+(\lambda_0,\ldots,\lambda_{n-1})$ give maps
$$
\op{Hom} ( \lambda_0,\lambda_1) \otimes \ldots \otimes \op{Hom}
(\lambda_{n-2},\lambda_{n-1} ) \to \op{Hom}
(\lambda_0,\lambda_{n-1})
$$
which are of degree $n-3$, when $n \ge 3$. These correspond to the $A_\infty$
multiplications $m_{n-1}$ when $n \ge 3$  (when appropriate sign conventions,
and orientations on the cells corresponding to
$D^+(\lambda_0,\ldots,\lambda_{n-1})$ are chosen).  The formula for the
differentials $ \d D$ gives the $A_\infty$ relation; indeed this is essentially
the original definition of $A_\infty$ algebra of Stasheff \cite{sta1963}.

When $n = 2,1$ the maps $D^+(\lambda)$ and $D^+(\lambda,\lambda')$
are of degree $0$.  $D^+(\lambda)$ gives a map $\K \to
\op{Hom}(\lambda,\lambda)$, which gives the unit in the $A_\infty$
category.   The axioms for units in an $A_\infty$ category
correspond to the relations in $\choc^+_{open,\Lambda}$ described
in \ref{lemma generators relations operad}.

$D^+(\lambda,\lambda')$ is the identity map $\op{Hom} (\lambda,
\lambda') \to \op{Hom}( \lambda,\lambda')$.
\end{proof}

\begin{lemma}
A split monoidal functor $\Phi : \choc^d_{open,\Lambda} \to \Comp$ is  the same
as a unital Calabi-Yau $A_\infty$ category with set of objects $\Lambda$.
\end{lemma}
\begin{proof}
This follows from the generators and relations description for the
category $\choc^d_{open,\Lambda}$.  There are two more generators
for $\choc^d_{open,\Lambda}$ over $\choc^+_{open,\Lambda}$, namely
the disc with two incoming and two outgoing boundaries.  These
give the pairing $\op{Hom}(\lambda_0,\lambda_1) \otimes \op{Hom}(
\lambda_1,\lambda_0) \to \K[d]$, and its inverse.  The extra
relations in $\choc^d_{open}$ correspond to the cyclic symmetry
condition.

\end{proof}

\begin{definition}
A unital extended Calabi-Yau $A_\infty$ category, with objects
$\Lambda$, is a h-split symmetric monoidal functor $\Phi :
\choc^d_{open,\Lambda} \to \Comp$.
\end{definition}
So a split extended Calabi-Yau $A_\infty$ category is the same as
an ordinary Calabi-Yau $A_\infty$ category.  Let $\Phi$ be an
extended Calabi-Yau $A_\infty$ category. For each   $\alpha \in
\op{Ob} \choc^d_{open,\Lambda}$ there is a complex $\Phi(\alpha)$,
and quasi-isomorphisms $\Phi(\alpha ) \otimes \Phi(\beta) \to
\Phi( \alpha \amalg \beta)$.  There are maps
$$
\Phi ( \alpha \amalg \{ \lambda_0,\ldots,\lambda_{n}\} ) \to \Phi
( \alpha \amalg \{\lambda_0,\lambda_{n}\})
$$
coming from the disjoint union of the disc
$D^+(\lambda_0,\ldots,\lambda_{n-1})$ and the identity map $\alpha
\to \alpha$.  These play the role of the $A_\infty$ operations
$m_{n}$, when $n \ge 2$. They satisfy relations analogous to the
usual $A_\infty$ relation. There are also maps
\begin{align*}
\Phi (\alpha \amalg \{\lambda_0,\lambda_1\}) & \to \Phi ( \alpha)
\\ \Phi ( \alpha ) & \to \Phi ( \alpha \amalg
\{\lambda_0,\lambda_1\})
\end{align*}
which play the role of the pairing and its inverse.  A cyclic
symmetry condition holds for the operation $\Phi ( \alpha \amalg
\{ \lambda_0,\ldots,\lambda_{n-1},\lambda_0 ) \to \Phi(\alpha)$
constructed from the $A_\infty$ operation $m_n$ and the pairing.
Also there are units, in $\Phi ( \{\lambda,\lambda\})$ satisfying
the usual constraints.

\begin{lemma}
The category of unital extended Calabi-Yau $A_\infty$ categories of
dimension $d$, with set of objects $\Lambda$, is quasi-equivalent,
in the sense of definition \ref{def quasi_equiv} to the category of
open TCFTs of dimension $d$.
\end{lemma}
This is immediate from theorem \ref{theorem quasi equivalence
module}, and the fact that $\choc^d_{open,\Lambda}$ is
quasi-isomorphic to $\mscr{O}^d_\Lambda$.  Thus we have proved
theorem A, part 1.

\begin{definition}
A unital extended $A_\infty$ category, with set of objects
$\Lambda$, is a h-split monoidal functor $\Phi :
\choc^+_{open,\Lambda} \to \Comp$.
\end{definition}
This makes sense, as we have already seen that such a functor
which is split is the same as a unital $A_\infty$ category.
\begin{proposition}
The following categories are quasi-equivalent.
\begin{enumerate}
\item
The category of unital extended $A_\infty$ categories, with set of
objects $\Lambda$.
\item
The category of unital  $A_\infty$ categories, with set of objects
$\Lambda$.
\item
The category of unital dg categories, with set of objects
$\Lambda$.
\end{enumerate}
\end{proposition}
\begin{proof}

Observe that for each pair of objects $\alpha,\beta \in \op{Ob}
\choc^+_{open,\Lambda}$, $H_i ( \choc^+_{open,\Lambda}(\alpha,\beta)) = 0$ if
$i \neq 0$.  As, the morphisms spaces are chains on moduli spaces of marked
points on discs, which are contractible. Also the complexes
$\choc^+_{open,\Lambda}(\alpha,\beta)$ are concentrated in degrees $\ge 0$.
This implies that $\choc^+_{open,\Lambda}$ is formal, that is quasi-isomorphic
to its homology, and quasi-isomorphic to $H_0 ( \choc^+_{open,\Lambda})$.

It is not difficult to see that a split functor $\Phi : H_0 (
\choc^+_{open,\Lambda} ) \to \Comp$ is the same as a unital dg
category with set of objects $\Lambda$.  Indeed, $ H_0 (
\choc^+_{open,\Lambda} (\{\lambda_0,\lambda_1,\ldots,\lambda_n\},
\{\lambda_0,\lambda_n\})$ is one dimensional, and corresponds to
the product map
$$
\op{Hom} ( \lambda_0,\lambda_1) \otimes \ldots \otimes \op{Hom}
(\lambda_{n-1},\lambda_n) \to \op{Hom} (\lambda_0,\lambda_n)
$$
Because it is one-dimensional, associativity holds.  Further, $H_0
( \choc^+_{open,\Lambda}(\alpha,\beta))$ is given by disjoint
unions of morphisms of this type.

Call a h-split functor $\Phi : H_0 ( \choc^+_{open,\Lambda} ) \to
\Comp$ a unital extended category.  Since there is a
quasi-isomorphism
$$
 \choc^+_{open,\Lambda}  \to H_0 ( \choc^+_{open,\Lambda} )
$$
there is a quasi-equivalence between unital extended $A_\infty$
categories and unital extended categories.

It remains to show how to remove the adjective extended. The
category $ \choc^+_{open,\Lambda} $ has the property that the maps
\begin{multline*}
\choc^+_{open,\Lambda} (\alpha_1,\{\lambda_1,\lambda_1'\}) \otimes
\choc^+_{open,\Lambda} (\alpha_2,\{\lambda_2,\lambda_2'\}) \ldots
\otimes \choc^+_{open,\Lambda} (\alpha_n,\{\lambda_n,\lambda_n'\})
\\ \to \choc^+_{open,\Lambda} (\alpha_1 \amalg \ldots \amalg \alpha_n
,\{\lambda_1,\lambda_1'\} \amalg \ldots \amalg
\{\lambda_1,\lambda_1'\}) \tag{\dag} \label{eqn split category}
\end{multline*}
 are isomorphisms.

Let $\Phi$ be a unital extended $A_\infty$ category.  Define a
unital $A_\infty$ category $F(\phi)$, i.e.\ a split monoidal
functor $\choc^+_{open,\Lambda} \to \Comp$, by
$$
F(\phi) ( O,s,t) = \otimes_{i = 0}^{O - 1}  \Phi ( \{s(i),t(i)\})
$$
There are maps $F(\phi) (\alpha) \to \phi(\alpha)$. Composing with the action
of $\choc^+_{open,\Lambda}$ on $\phi$ gives maps
$$
F(\phi) (\alpha) \otimes \choc^+_{open,\Lambda}
(\alpha,\{\lambda_0,\lambda_1\}) \to F(\phi) (
\{\lambda_0,\lambda_1\})
$$
Because the maps \ref{eqn split category} are isomorphisms, it
follows that these extend to give a unique $\choc^+_{open,\Lambda}
$ module structure on $F(\phi)$, with a quasi-isomorphism $F(\phi)
\to \phi$.

This shows that the category of extended $A_\infty$ categories is
quasi-equivalent to the category of $A_\infty$ categories.
Similarly the category of extended dg categories is
quasi-equivalent to the category of dg categories; this finishes
the proof.

\end{proof}

The proof shows something stronger; the obvious map from dg
categories to $A_\infty$ categories is half of a
quasi-equivalence.  This means that every $A_\infty$ category is
quasi-isomorphic, in a functorial way, to a dg category.

\subsection{The Hochschild chain complex}
For an associative algebra $A$, over our ground field $\K$, and an
$A$-bimodule $M$, recall the Hochschild complex $C_\ast(A,M)$ is
defined by
$$
C_n(A,M) = M \otimes A^{\otimes n}
$$
The differential $\d : C_n(A,M) \to C_{n-1}(A,M)$ is given by the
formula
\begin{multline} \d ( m \otimes a_1 \otimes \ldots  a_n )
= m a_1 \otimes a_2 \otimes \ldots a_n  \\ + \sum_{i = 1}^{n-1}
(-1)^i m \otimes a_1 \otimes \ldots a_{i} a_{i+1} \otimes \ldots
a_n + (-1)^n a_n m \otimes a_1 \otimes \ldots a_{n-1}
\end{multline}
When $M = A$, we write $C_\ast(A)$ for $C_\ast(A,A)$.

The normalised Hochschild chain complex is a quotient of
$C_\ast(A,M)$ by the contractible complex spanned by elements $m
\otimes a_1 \otimes \ldots a_n$ where at least one of the $a_i =
1$.   We write $\br C_\ast(A,M)$ for the normalised chain complex,
and $\br C_\ast(A)$ for the normalised chain complex with
coefficients in $A$.

Similar definitions hold for dg algebras $A$ and dg modules $M$,
except with  extra terms in the differential coming from the
differential on $A$ and $M$, and a change in sign coming from the
grading on $A$ and $M$.

Let $A$ be a dg category.  Define the Hochschild chain complex
$$
C_\ast(A) = \oplus  ( \op{Hom} ( \alpha_0,\alpha_1) \otimes \ldots
\otimes \op{Hom} ( \alpha_{n-1},\alpha_0) ) [1-n])
$$
where the direct sum is over $n$ and sequences $\alpha_0,\ldots,
\alpha_{n-1}$ of objects in $A$.

The differential is given by essentially the same formula as in
the algebra case:
\begin{multline*}
\d (\phi_0 \otimes \ldots \otimes \phi_{n-1}) = \sum_{i = 0}^{n-1} \pm \phi_0 \ldots \d \phi_i \ldots \phi_{n-1} \\
+  \sum_{i = 0}^{n-2} \pm  \phi_0 \ldots \otimes  (\phi_{i+1} \circ \phi_{i})
\otimes \ldots \otimes \phi_{n-1} \pm(\phi_{0}\circ \phi_{n-1}) \otimes \ldots
\otimes \phi_{n-2}
\end{multline*}
If $A$ is unital, then we can define the normalised Hochschild
 chain complex $\br C_\ast(A)$ by taking the quotient by the
contractible subcomplex spanned by $\phi_0 \otimes \ldots
\phi_{n-1}$ where at least one of the $\phi_i$, where $i > 0$, is
an identity map.

\begin{lemma}
The functor $A \mapsto \br C_\ast(A)$ is an exact functor from the
category of dg categories with fixed set of objects $\Lambda$ to
the category of complexes.
\end{lemma}
\begin{proof}
Give the normalised Hochschild chain complex $\br C_\ast(A)$ the
obvious filtration, defined by $F^i ( \br C_\ast(A) )$ is the
subcomplex spanned by $\phi_0 \otimes \ldots \phi_{i-1}$.   If $A
\to B$ is a map of dg categories with fixed set of objects, the
induced map $\br C_\ast(A) \to \br C_\ast(B)$ preserves the
filtration.  We need to show that if $A \to B$ is a
quasi-isomorphism then so is $\br C_\ast(A) \to \br C_\ast(B)$. It
is sufficient to show that the associated graded map is a
quasi-isomorphism; but this is obvious.
\end{proof}

\begin{definition}
Let $A$ be a (possibly extended) $A_\infty$ category.  Define the
Hochschild homology $HH_\ast(A)$ to be the homology of the dg
category associated to it under the quasi-equivalence between
(extended) $A_\infty$ and dg categories.

If $\Phi$ is an extended Calabi-Yau $A_\infty$ category, define
the Hochschild homology of $\Phi$  to be the homology of the
associated extended $A_\infty$ category.
\end{definition}
We could also use an explicit complex to define the Hochschild
homology, but this would involve getting the signs correct.

\begin{proposition}
Let $\Phi$ be a unital extended Calabi-Yau $A_\infty$ category.
Then
$$
H_\ast ( \choc^d(-,1)_{\Lambda} \otimes_{\choc^d_{open,\Lambda}}
\Phi ) = HH_\ast(\Phi)
$$
\end{proposition}
\begin{proof}

Recall the definition of $\choc^+$ in definition \ref{definition
dplus}. We have a generators and relations description of
$\choc^d$, in theorem \ref{theorem open closed generators
relations}, and we defined $\choc^+$ to have the same generators
and relations but as a $\op{Ob} \mscr{OC}-\choc^+_{open}$ bimodule
rather than a $ \op{Ob} \mscr{OC}-\choc^d_{open}$ bimodule.

We have
$$
\choc^d(-,1)_{\Lambda} \otimes_{\choc^d_{open,\Lambda}} \Phi =
\choc^+(-,1)_{\Lambda} \otimes_{\choc^+_{open,\Lambda}} \Phi
$$
Further, we have shown that the functor $\choc^+(-,1))_{\Lambda}
\otimes_{\choc^+_{open,\Lambda}}  -$ is exact (lemma \ref{lemma
flat and split}).  The $\choc^+_{open,\Lambda}$ module underlying
$\Phi$ is the extended $A_\infty$ category associated to $\Phi$.

What remains to be shown is that, for an actual dg category $B$,
considered as a left $\choc^+_{open,\Lambda}$ module,
$$
H_\ast ( \choc^+(-,1)_{\Lambda} \otimes_{\choc^+_{open,\Lambda}} B ) =
HH_\ast(B)
$$
We will show something a bit more; we will show that
$$
\choc^+(-,1)_{\Lambda} \otimes_{\choc^+_{open,\Lambda}} B = \br
C_\ast(B)
$$
is the normalised Hochschild chain complex.

This follows from the generators and relations description of the
right $\choc^+_{open,\Lambda}$ module, $\choc^+(-,1)_{\Lambda}$.
Recall it is generated by the annuli
$A(\lambda_0,\ldots,\lambda_{n-1})$, modulo the relation that when
we glue the disc with one outgoing marked point onto a marked
point of $A(\lambda_0,\ldots,\lambda_{n-1})$, we get zero, except
for the marked point between $\lambda_{n-1}$ and $\lambda_0$. The
annulus $A(\lambda_0,\ldots,\lambda_{n-1})$ is in degree $n-1$.

This shows us that  $\choc^+(-,1)_{\Lambda}
\otimes_{\choc^+_{open,\Lambda}} B$, as a vector space, is the
quotient of $\oplus  ( \op{Hom} ( \alpha_0,\alpha_1) \otimes
\ldots \otimes \op{Hom} ( \alpha_{n-1},\alpha_0) ) [1-n])$ by the
subspace spanned by elements of the form $\phi_0 \otimes \ldots
\phi_{n-1}$, where at least one of the $\phi_i$ with $i > 0$ is an
identity map.

That is, as a vector space, there is a natural isomorphism
$$
\choc^+(-,1)_{\Lambda} \otimes_{\choc^+_{open,\Lambda}} B \iso \br
C_\ast (B)
$$
It remains to show that this is compatible with the differential.
This follows immediately from the formula for the differential of
the annulus $A(\lambda_0,\ldots,\lambda_{n-1})$, see figure
\ref{figure differential annuli}.  Recall that $m_n = 0$ when $n >
2$ in our category $B$, so that the only thing that contributes is
when two marked points on the boundary of the annulus collide.  This
corresponds to composing the corresponding consecutive morphisms in
the formula for the Hochschild differential.

\end{proof}

This completes the proof of theorem A.

We have shown theorem A part 1: the category of unital extended
Calabi-Yau $A_\infty$ categories is quasi-equivalent to the
category of open TCFTs.

We have also shown that under the quasi-equivalence between
$\op{Ob} \mscr{OC}^d_{\Lambda} - \choc^d_{open,\Lambda}$ bimodules
and $\op{Ob} \mscr{OC}^d_{\Lambda} - \mscr{O}^d_{\Lambda}$
bimodules, the bimodule $\choc^d_{\Lambda}$ corresponds to
$\mscr{OC}^d_{\Lambda}$.  Also, $\choc^d_{\Lambda}$ is flat.

Thus, by lemma \ref{lemma tensor compatible} if $M$ is a left
$\choc^d_{open,\Lambda}$ module corresponding to a left
$\mscr{O}^d_{\Lambda}$ module $M'$, then
$$
\mscr{OC}^d_{\Lambda} (-,\beta) \otimes^{\mbb
L}_{\mscr{O}^d_{\Lambda}} M' \iso \choc^d_{\Lambda}( - ,\beta)
\otimes_{\choc^d_{open,\Lambda}} M
$$
Denote by $N(\beta)$ the left hand side of this equation. Then
$N(\beta)$ is h-split, if $M$ is; the maps $N(\beta) \otimes
N(\beta') \to N(\beta \amalg \beta')$ are quasi-isomorphisms. This
shows that $N$ defines an open-closed TCFT of dimension $d$, which
is the homotopy universal open-closed TCFT associated to $M'$.

Finally, we have calculated the homology of the closed states of
$N$ to be the Hochschild homology of the associated $A_\infty$
category.

\section{Appendix}
In this appendix a symmetric monoidal functor $C_\ast$ from the category of topological
spaces with local systems to chain complexes is constructed, which computes
homology groups, and satisfies several nice properties. In particular, for a
cell complex $X$, there is a map $C_\ast^{cell}(X) \to C_\ast(X)$ which is
natural for a strong notion of cellular map.

We recall the properties of homology with local coefficients. A
$\K$ local system on a space $Y$ is a locally constant sheaf of
$\K$ vector spaces on $Y$. If $E$ is a local system on $Y$, there
are homology groups $H_i(Y,E)$ with local coefficients.   Spaces
with local systems form a category; a map $(Y,E) \to (Z,F)$ is a
map $f : Y \to Z$ and a map $E \to f^\ast F$.  Homology with
coefficients defines a functor from this category to the category
of graded $\K$ vector spaces.

This functor satisfies the following properties.
\begin{enumerate}
\item
If
$$
0 \to E_1 \to E_2 \to E_3 \to 0
$$
is a short exact sequence of local systems on $Y$, there is a corresponding
long exact sequence of homology groups
$$
\ldots \to H_i(E_1) \to H_i(E_2) \to H_i(E_3) \to H_{i-1} (E_1) \to \ldots
$$
\item
If $Y = U \cup V$ is written as a union of open subsets, there is a
Mayer-Vietoris long exact sequence of homology groups
$$
\ldots \to H_i(U \cap V, E) \to H_i(U,E) \oplus H_i(V,E) \to H_i(Y,E) \to
H_{i-1}(U \cap V, E) \to \ldots
$$
\item
Two maps $f_0,f_1 : (Y,E) \to (Z,F)$ are homotopic if they extend to a map $F :
(Y\times I, \pi_1^\ast E) \to (Z,F)$.  Homotopic maps induce the same map on
homology.
\item
If $Y = \ast$ is a point, and $E$ is a vector space, then $H_i(\ast, E) = 0$ if
$i  \neq 0$, and $H_0(\ast,E) = E$.
\end{enumerate}
On reasonable spaces, for example spaces with the homotopy type of
finite cell complexes, the functor $(Y,E) \to H_\ast(Y,E)$ is
determined by these properties.  We can define $H_i(Y,E)$ using
singular simplices $f : \diag_n \to Y$ with sections of $f^\ast E
\otimes \omega$, where $\omega$ is the orientation sheaf. There
are also relative homology groups $H_i(Y,Y',E)$ for a subspace $Y'
\subset Y$ and a local system $E$ on $Y$, which fit into the
obvious long exact sequence.

A finite regular cell complex is a space $X$ obtained by attaching finitely
many cells to a finite number of points, with the property that the boundary of
one cell is a union of lower dimensional cells.  Let $X_i \subset X$ be the
union of cells of dimension $\le i$. A strong cellular map between finite
regular cell complexes $X, X'$ is a continuous  map $f : X \into X'$ such that
$f^{-1}(X'_i) = X_i$.    Thus we have a category $\op{Cell}$ of finite regular
cell complexes with these morphisms.

For a topological space $Y$, let $\op{Cell}_Y$ be the category whose objects
are finite regular cell complexes $X$ with a map $f : X \to Y$, and whose
morphisms are strong cellular maps $X \to X'$ such that the obvious diagram
commutes.

There is a functor $C_\ast^{cell} : \op{Cell}_Y \to \Comp$, which
takes $X$ to the $\K$ cellular chain complex
$C_\ast^{cell}(X,\K)$. (Of course we could use any coefficient
ring).  If $E$ is a local system on $Y$, then pulling back $E$
gives a local system on each object $X \in \op{Cell}_Y$, and there
is a functor of cellular chains with coefficients from
$\op{Cell}_Y \to \Comp$. This functor applied to $X \in
\op{Cell}_Y$ is denoted $C_\ast^{cell}(X,E)$. By definition,
$C_n^{cell}(X,E) = H_n(X_n,X_{n-1}, E)$ is the relative sheaf
homology. $C_n^{cell}(X,E)$ is naturally isomorphic to the space
of sections over $X_n \setminus X_{n-1}$ of the sheaf $E \otimes
\omega_{X_n \setminus X_{n-1}}$, where $\omega_{X_n\setminus
X_{n-1}}$ is the orientation sheaf.

Define $C_\ast(Y, E)$ by
$$
C_\ast(Y,E) = \lim_{\substack{\to \\ X  \in \op{Cell}_Y }} C_\ast^{cell}(X,E)
$$
to be the direct limit over the cellular chain groups of objects of
$\op{Cell}_Y$.

It is clear that $C_\ast$ is functorial.  Denote by $H'(Y,E)$ the homology of
the chain complex $C_\ast(Y,E)$.
\begin{proposition}
The functor $H'(Y,E)$ satisfies the axioms (1)-(4) above, and so coincides with
usual homology with local coefficients on reasonable spaces.
\end{proposition}
\begin{proof}
Axiom (1) is straightforward; the sequence of complexes
$$
0 \to C_\ast(Y,E_1) \to C_\ast(Y,E_2) \to C_\ast(Y,E_3) \to 0
$$
is exact. Axiom (3) follows as if $\iota_0,\iota_1 : Y \into Y
\times I$ are the inclusions, there is a canonical chain homotopy
between the induced maps $C_\ast(Y,E) \to C_\ast(Y \times
I,\pi_1^\ast E)$.  Axiom (4) is also quite straightforward; for
any $n$ cell complex $X$ over a point, with $n > 1$, there is a
cellular isomorphism $X \to X$ changing the orientation on the $n$
cells.

It remains to prove the Mayer-Vietoris axiom.  The sequence of complexes
$$
0 \to C_\ast(U \cap V,E) \to C_\ast(U,E) \oplus C_\ast(V,E) \to C_\ast(Y,E) \to
0
$$
is actually exact.  Exactness on the left and in the middle is straightforward.
Exactness on the right is more difficult; this can be proved by showing,
inductively on the dimension of the cells,  that for any $n$ cell complex $X
\in \op{Cell}_Y$, we can find a refinement $X'$ of the cell structure on $X$
such that any closed cell of $X'$ lands in either $U$ or $V$.

\end{proof}

If $X_1,X_2$ are cell complexes, and $E_i$ are finite dimensional
$\K$ local systems on $X_i$, then there is an isomorphism
$$
C_\ast^{cell}(X_1 , E_1) \otimes C_\ast^{cell}(X_2,E_2) \xrightarrow{\iso}
C_\ast^{cell}( X_1 \times X_2, E_1 \boxtimes E_2)
$$
This induces maps
$$C_\ast(Y_1,
E_1) \otimes C_\ast(Y_2,E_2) \to C_\ast(Y_1 \times Y_2, E_1
\boxtimes E_2)$$ making $C_\ast$ into a symmetric monoidal functor from spaces
with finite dimensional $\K$ local systems to chain complexes.

\def\cprime{$'$}

\end{document}